\documentclass[11pt]{article}
\usepackage{amsfonts,latexsym,amsmath,amscd,geometry}
\usepackage{amssymb}
\usepackage{latexsym}
\usepackage{color}

\newcommand \nc{\newcommand}
\newtheorem{theorem}{Theorem}[section]
\newtheorem{lemma}[theorem]{Lemma}
\newtheorem{proposition}[theorem]{Proposition}
\newtheorem{corollary}[theorem]{Corollary}

\newtheorem{remark}[theorem]{Remark}

\nc{\ba}{\begin{array}}\nc{\ea}{\end{array}}
\nc{\be}{\begin{eqnarray}}\nc{\ee}{\end{eqnarray}}
\nc{\beq}{\begin{equation}}\nc{\eeq}{\end{equation}}
\nc{\bex}{\begin{eqnarray*}}\nc{\eex}{\end{eqnarray*}}
\nc{\btm}{\begin{theorem}} \nc{\etm}{\end{theorem}}
\nc{\blm}{\begin{lemma}} \nc{\elm}{\end{lemma}}
\nc{\R}{\mathbb{R}}  \nc{\ld}{\lambda}
\nc{\va}{\varphi}
\nc{\ve}{\varepsilon}

\def\pf{\noindent{\bf Proof.\quad}}

\newcommand \qed {\hfill $\Box$}

\begin{document}
\title{Higher dimensional Ginzburg-Landau equations under weak anchoring boundary conditions
\author{Patricia Bauman, \ Daniel Phillips, \ Changyou Wang\\
Department of Mathematics\\
Purdue University\\
West Lafayette, IN 47907, USA}}
\maketitle
\begin{abstract} For $n\ge 3$ and  $0<\epsilon\le 1$,  let $\Omega\subset\mathbb R^n$ be a bounded smooth domain and $u_\epsilon:\Omega \subset\R^n\to \mathbb R^2$
solve the Ginzburg-Landau equation under the weak anchoring boundary condition:
$$\begin{cases}
-\Delta u_\epsilon=\frac{1}{\epsilon^2}(1-|u_\epsilon|^2)u_\epsilon &\ {\rm{in}}\ \ \Omega,\\
\frac{\partial u_\epsilon}{\partial\nu}+\lambda_\epsilon(u_\epsilon-g_\epsilon)=0 & \ {\rm{on}}\ \ \partial\Omega,
\end{cases}
$$
where the anchoring strength parameter $\lambda_\epsilon=K\epsilon^{-\alpha}$ for some $K>0$ and $\alpha\in [0,1)$,
and $g_\epsilon\in C^2(\partial\Omega, \mathbb S^1)$. Motivated by the connection with the Landau-De Gennes model
of nematic liquid crystals under weak anchoring conditions,  we study the {asymptotic behavior} of $u_\epsilon$  as
$\epsilon$ goes to zero under the condition that the total modified Ginzburg-Landau energy {satisfies} $F_\epsilon(u_\epsilon,\Omega)\le M|\log\epsilon|$
for some $M>0$.
\end{abstract}

\section{Introduction}
\setcounter{equation}{0}
\setcounter{theorem}{0}

Given a bounded smooth domain $\Omega\subset\R^n$ with $n\ge 2$, the Ginzburg-Landau energy
for a map $u:\Omega\to \mathbb R^2$ is defined by
$$E_\epsilon(u,\Omega)=\int_\Omega e_\epsilon(u), \ {\rm{with}}\ e_\epsilon(u)\equiv
\frac12|\nabla u|^2+\frac{1}{4\epsilon^2}(1-|u|^2)^2 \ {\rm{and}}\ 0<\epsilon\le 1.$$
Since the pioneering work by Brezis-Bethuel-Helein
\cite{BBH, BBH1} in dimension $n=2$, there have been extensive studies {on the asymptotic behavior} of minimizers or critical points $u_\epsilon$
of the Ginzburg-Landau energy $E_\epsilon$ in $\Omega$, under the Dirichlet boundary condition $g_\epsilon$,
as $\epsilon$ goes to zero. Note that {any such $u_\epsilon$} is a smooth
solution of the Ginzburg-Landau equation under the Dirichlet boundary condition:
\begin{equation}\label{GL-Dir}\displaystyle
\begin{cases}
-\Delta u_\epsilon=\frac{1}{\epsilon^2}(1-|u_\epsilon|^2)u_\epsilon & \ {\rm{in}}\  \Omega,\\
u_\epsilon= g_\epsilon & \ {\rm{on}}\ \partial\Omega.
\end{cases}
\end{equation}
Among other results, it was shown {in \cite{BBH, BBH1} that for energy minimizers when $n=2$,} 
if $g_\epsilon=g\in C^\infty(\partial\Omega,\mathbb S^1)$ and ${\rm{deg}}(g)=d$, then {there exists}
a vortex set $\Sigma_*\subset\Omega$ of exactly $|d|$ points and a smooth harmonic map $u_*\in C^\infty(\Omega\setminus\Sigma_*, \mathbb S^1)$, with $u_*=g$ on $\partial\Omega$, such that after taking a subsequence, $u_\epsilon\rightarrow u_*$ in $C^\infty_{\rm{loc}}(\Omega\setminus\Sigma_*,\mathbb R^2)\cap W^{1,p}(\Omega,\mathbb R^2)$ {for any $1\le p<2$. (See also Struwe \cite{Struwe1, Struwe2}.)} A similar result holds for general solutions $u_\epsilon$ of \eqref{GL-Dir}, with $\Sigma_*$ having finite but not necessarily $|d|$ points, provided $E_\epsilon(u_\epsilon, \Omega)\le C|\log\epsilon|$.  When $n\ge 3$,
the {asymptotic behavior of minimizing solutions $u_\epsilon$ to \eqref{GL-Dir} has been} studied by
Rivi\`ere \cite{Riviere} ($n=3$), Lin-Rivi\`ere \cite{LR1, LR1.0}, Sandier \cite{Sandier} ($n=3$),
Alberti-Baldo-Orlandi \cite{ABO} and Jerrard-Soner \cite{JS}. It was shown that the vortex set $\Sigma_*$ is ($n-2$)-dimensional area minimizing,
with $\partial\Sigma_*\subset\partial\Omega$, and $u_\epsilon\rightarrow u_*$ in $C^\infty_{\rm{loc}}(\Omega\setminus\Sigma_*,\mathbb R^2)\cap W^{1,p}(\Omega, \mathbb R^2)$ for $1\le p<\frac{n}{n-1}$. {The asymptotic behavior of non-minimizing
solutions to \eqref{GL-Dir} has} been studied by Lin-Rivi\`ere \cite{LR2} for $n=3$ and Brezis-Bethuel-Orlandi \cite{BBO1, BBO2} for $n\ge 3$,
in which the vortex set $\Sigma_*\subset\overline\Omega$ is shown to be a ($n-2$)-rectifiable set.
{We would also like to mention the extensive studies for the Ginburg-Landau equation with magnetic fields by Sandier-Serfaty \cite{SandierSerfaty}.}

The Ginzburg-Landau equation under the weak anchoring boundary condition arises, if we impose the boundary behavior
through the addition of a surface energy term into $E_\epsilon$:
\begin{equation}
F_\epsilon(u, \Omega)=\int_\Omega \frac12|\nabla u|^2+\frac{1}{4\epsilon^2}(1-|u|^2)^2
+\frac{\lambda_\epsilon}2\int_{\partial\Omega}|u-g_\epsilon|^2,
\end{equation}
{where $\lambda_\epsilon>0$ is called an} anchoring strength parameter. It is readily seen that the Euler-Lagrange equation for
a critical point $u_\epsilon$ of $F_\epsilon$ is
\begin{equation}\label{GL-WA0}
\begin{cases}
-\Delta u_\epsilon=\frac{1}{\epsilon^2}(1-|u_\epsilon|^2)u_\epsilon & \ {\rm{in}}\ \Omega,\\
\frac{\partial u_\epsilon}{\partial\nu}+\lambda_\epsilon(u_\epsilon-g_\epsilon)=0 & \ {\rm{on}}\ \partial\Omega.
\end{cases}
\end{equation}
Here $\nu$ denotes the unit outward normal of $\partial\Omega$.

This type of boundary value problem is closely related to interesting problems arising
from the study of bulk nematic liquid crystals with oil droplets or nano-particles included
(see, for example, Kleman-Lavrentovich \cite{KL}), where the boundary behavior at the interface
between the droplet and the bulk nematic liquid crystals is usually constrained by surface energies rather than
a prescribed Dirichlet data (called the strong anchoring condition). Recall that one of the
most universal models to describe nematic liquid crystals is the Landau-De Gennes model
\cite{De-P, MN}, in which De Gennes proposed to represent non-oriented direction fields of liquid crystals
by symmetric traceless $n\times n$ matrix-valued functions $Q(x)$, called $Q$-tensors.
Note that the two widely used simplified models -- the Oseen-Frank model \cite{HKL} and the Ericksen model \cite{Ericksen},
which utilize unit vector fields $d:\Omega\to\mathbb S^{n-1}$ to describe nematic liquid crystals, can
be embedded into the Landau-De Gennes model via the identification $Q(x)=s(d\otimes d-\frac{1}{n} \mathbb I_n)$,
$s\in\mathbb R$, called the uniaxial $Q$-tensor. {Recall that a simplified} version
of the Landau-De Gennes functional, with weak anchoring boundary conditions, takes the form \cite{KL, MN}:
$$E_{\rm{LdG}}(Q)=\int_\Omega \big(\frac12|\nabla Q|^2+\frac{1}{L}f_b(Q)\big)\,dx+W\int_{\partial\Omega}\frac12|Q-Q^0|^2\,dH^{n-1},$$
where $L>0$ is the elasticity constant, $W>0$ is the relative anchoring strength constant, $Q^0$ is a prescribed $Q$-tensor
function preferred by liquid crystal materials on $\partial\Omega$, and the bulk potential function
$$f_b(Q)=-\frac{a}2 {\rm{tr}}(Q^2)-\frac{b}3 {\rm{tr}}(Q^3)+\frac{c}4 {\rm{tr}}^2(Q^2)-d,$$
which penalizes $Q$ for not being uniaxial. In fact, if $a, b, c>0$ and $d$ is chosen so that
min $f_b=0$, then $f_b$ is minimal iff
$$Q=s_+(d\otimes d-\frac{1}{n}\mathbb I_n),$$
for some specific scalar parameter constant $s_+=s_+(a,b,c)>0$. Similar to the fact that the Ginzburg-Landau energy $E_\epsilon$ approximates the Dirichlet energy of harmonic maps to $\mathbb S^1$, the Landau-De Gennes functional is a relaxation of the Dirichlet energy of harmonic maps  to the space of uniaxial $Q$-tensor fields. {For a planar sample ($n=2$) where the director field lies in the same plane of the sample, Majumdar \cite{Majumdar} and Ball-Zarnescu \cite{BZ1, BZ2} showed equivalence between the Landau-De Gennes  model for $2\times2$ $Q$-tensor fields and the Ginzburg-Landau model for
complex-valued functions.
In particular they showed that in dimension $n=2$,
any vortex of integer degree $k$ in a solution to the Ginzburg-Landau energy corresponds to a vortex of degree $\frac{k}2$ in a solution to the $Q$-tensor functionals.
Such a vortex  can also be viewed as the cross-section of a disclination line singularity in dimension three.
If the director field is not constrained to be planar, the issue that is studied is the asymptotics of minimizers $Q_L$ of the Landau-De Gennes functional $E_{\rm{LdG}}$ among functions valued in the space of $3\times3$ trace-free matrices,
as $L\to 0$.  In \cite{BPP} Bauman, Park, and Phillips analyzed minimizers for this energy
in dimension $n=2$ among $3\times3$ trace-free $Q$-tensor fields in which $\vec{e}_3$ is an eigenvector, which models thin films; they showed that for $L$ sufficiently small and appropriate Dirichlet data, all vortices of minimizers have degree $\frac{1}2$.
The asymptotics of minimizers $Q_L$ of $E_{\rm{LdG}}$ among $3\times3$ symmetric trace-free tensors
as $L\to 0$ was investigated (under Dirichlet boundary conditions) in dimension $n=2$ by Golovaty-Montero \cite{GM}, Canevari \cite{C1} and in $n=3$ by Majumdar-Zarnescu \cite{MZ} and Canevari \cite{C2}.  In these studies,
defects of the ``hedgehog" type and some evidence of disclination line defects in certain settings were obtained.}

In an interesting article \cite{ABG}, Alama-Bronsard-Galav\~ao-Sosua {studied in dimension $n=2$}
the asymptotics of minimizers $u_\epsilon$ of the Ginzburg-Landau energy $F_\epsilon$ with weak anchoring
conditions, when the anchoring strength parameter takes a prescribed rate $\lambda_\epsilon=K\epsilon^{-\alpha}$
for some $K>0$ and $\alpha\in (0,1)$. In particular, they {demonstrated the effect} of weak anchoring conditions by
showing that the set $\Sigma_*$ of $|d|$-vortex points all lie inside $\Omega$ when $\alpha\in (\frac12,1)$; {they lie on} $\partial\Omega$
when $\alpha\in (0, \frac12)$; and there exists a theory of renormalized energy functions associated with the vortex points, similar to
that under the Dirichlet boundary condition by \cite{BBH} .

In this paper, we are mainly interested in deriving {asymptotics of general solutions} $u_\epsilon$ to the Ginzburg-Landau
equation with the weak anchoring condition (\ref{GL-WA0}) for $g_\epsilon$ satisfying certain conditions to be specified below,
under the assumption that $F_\epsilon(u_\epsilon,\Omega)\le M|\log \epsilon|$ (e.g., the same energy threshold as that of minimizers of $F_\epsilon$). The motivation for doing this is twofold:
{the first} is to extend earlier results on solutions {of Ginzburg-Landau equations with Dirichlet boundary conditions}, especially the works by
\cite{BBO1, BBO2} and \cite{ABG}, to higher dimensional Ginzburg-Landau equations under weak anchoring conditions; {the second is to gain further insight} on how to study the Landau-De Gennes model under weak anchoring conditions in dimension three.

Before stating our main result, we would like to specify throughout this paper {an assumption} on the weak anchoring data $\{g_\epsilon\}\subset C^2(\partial\Omega, \mathbb S^1)$,
$0<\epsilon\le 1$,  {called the} \\
Condition ({\bf G}): There exists $C_0>0$ independent of $\epsilon$ such that
\begin{equation}\label{bdry_bound2}
\|g_\epsilon\|_{C^2(\partial\Omega)}\le C_0.
\end{equation}
It follows from the condition ({\bf G}) that there exists $g_*\in C^2(\partial\Omega, \mathbb S^1)$, with $\|g_*\|_{C^2(\partial\Omega)}\le C_0$,
such that, after taking a subsequence,
$$g_\epsilon\rightarrow g_* \ {\rm{in}}\ C^1(\partial\Omega,\mathbb S^1).$$

From now on, we always assume dimensions $n\ge 3$. Our main result is stated as follows. 
\begin{theorem}\label{compactness1} For any $g_\epsilon\in C^2(\partial\Omega, \mathbb S^1)$ satisfying the {condition {\rm{(}}{\bf G}{\rm{)}} and $\lambda_\epsilon=K\epsilon^{-\alpha}$ for some $K>0$} and $\alpha\in [0,1)$, let
$u_\epsilon\in C^\infty(\overline\Omega,\mathbb R^2)$ be a solution of \eqref{GL-WA0}
for $0<\epsilon\le 1$. Assume that there exists $M>0$ such that
\begin{equation}\label{global_bound}
F_\epsilon(u_\epsilon,\Omega)\le M|\log\epsilon|, \ \forall\epsilon\in(0,1].
\end{equation}
Then for all $1\le p<\frac{n}{n-1}$  there exists $C_p>0$, depending on $p,M, K, \alpha, C_0$
and $\Omega$, such that
\begin{equation}\label{W1p-bound}
\big\|\nabla u_\epsilon\big\|_{L^p(\Omega)}\le C_p, \ \forall\epsilon\in(0,1].
\end{equation}
Hence there exists a subsequence $\epsilon_i\rightarrow 0$ and a map $u_*\in W^{1,p}(\Omega, \mathbb S^1)$
for all $1\le p<\frac{n}{n-1}$ such that the following statements hold:\\
\noindent {\rm{(a)}} $u_{\epsilon_i}\rightharpoonup u_*$ in $W^{1,p}(\Omega,\mathbb R^2)$
for all $1\le p<\frac{n}{n-1}$, and  $u_*$ is a generalized harmonic map to $\mathbb S^1$ in the sense that
${\rm{div}}(\nabla u_*\times u_*)=0 \ \ {\rm{in}}\ \ \mathcal D'(\Omega).$ Moreover,
\begin{itemize}
\item[{\rm{(a1)}}] $u_*=g_*$ on $\partial\Omega$ in the trace sense when $0<\alpha<1$, and
\item[{\rm{(a2)}}] $\big(\frac{\partial u_*}{\partial\nu}-Kg_*\big)\times u_{*}=0$ on $\partial\Omega$ in the distribution sense
(see \eqref{distribution} below) when $\alpha=0$.
\end{itemize}

\noindent{\rm{(b)}} there exists a nonnegative Radon measure $\mu_*$ in $\overline\Omega$ such that
$$\frac{1}{|\log\epsilon_i|}\Big(e_{\epsilon_i}(u_{\epsilon_i})\,dx+\frac{\lambda_{\epsilon_i}}2|u_{\epsilon_i}-g_{\epsilon_i}|^2\,dH^{n-1}
_{\partial\Omega}\Big)\rightharpoonup\mu_*,
\ {\rm{as}}\ i\rightarrow \infty,$$
as convergence of Radon measures in $\overline\Omega$. Set $\Sigma={\rm{supp}}(\mu_*)\subset\overline\Omega$. Then
$H^{n-2}(\Sigma)<\infty$, $u_*\in C^\infty(\Omega\setminus \Sigma,\mathbb S^1)$, and $u_{\epsilon_i}\rightarrow u_*$ in
$C^k_{\rm{loc}}(\Omega\setminus\Sigma,\mathbb R^2)$ for any $k\ge 1$;
$\Sigma$ is a $(n-2)$-rectifiable  set; and $\mu_*\ {\rm{L}}\ \Omega$ is a $(n-2)$-dimensional stationary varifold.
If, in addition, $\alpha=0$, then there exists $\delta\in (0,1)$ such that
$u_*\in C^\delta(\overline\Omega\setminus\Sigma, \mathbb S^1)$, and $u_{\epsilon_i}\rightarrow u_*$ in  $C^\delta_{\rm{loc}}(\overline\Omega\setminus\Sigma,\mathbb R^2)$.
\end{theorem}

Besides utilizing many previous techniques from \cite{BBO1, BBO2, LR2} that establish the interior estimates, such as the interior monotonicity
formula (\cite{BBO2} Lemma II.2) and the interior $\eta$-compactness\footnote{also called $\eta$-ellipticity by \cite{BBO2}} (\cite{BBO2} Theorem 2), one of the crucial ingredients that we need
in order to prove \eqref{W1p-bound} of Theorem \ref{compactness1} is to develop the $\eta$-compactness property
near $\partial\Omega$. More precisely, let $B_R(x)\subset\mathbb R^n$  be the ball with center $x\in\mathbb R^n$ and
radius $R>0$. Then
\begin{theorem}\label{bdry-eta-compact1} There exist $r_0=r_0(\Omega)>0$ and $\epsilon_0\in (0,1)$ depending on
$\eta, \Omega, K, \alpha$ and $C_0$ such that for any $\eta>0$,
$K>0$, {and $\alpha\in [0,1)$, if $\{g_\epsilon\} \subset C^2(\partial\Omega,\mathbb S^1)$ satisfies the condition {\rm{(}}{\bf G}{\rm{)}} and if}
$u_\epsilon\in C^2(\overline\Omega,\mathbb R^2)$ is a solution of \eqref{GL-WA0}
with $\lambda_\epsilon=K\epsilon^{-\alpha}$ that satisfies,
for any fixed $x_0\in\overline\Omega$, \begin{equation}\label{eta-condition}
\Phi_\epsilon(u_\epsilon, B_R(x_0)\cap\Omega):=R^{2-n}\Big(\int_{B_R(x_0)\cap\Omega}e_\epsilon(u_\epsilon)
+\int_{\partial\Omega\cap B_R(x_0)}\frac{\lambda_\epsilon}2|u_\epsilon-g_\epsilon|^2\Big)\le \eta|\log\epsilon|
\end{equation}
for $0<\epsilon\le\epsilon_0$ and $\epsilon^\alpha\le R\le r_0$,
then there exist $L>0$ and $\theta\in (0,1)$,  depending on $n, \Omega,
C_0, K, \alpha$, such that
\begin{equation}\label{bdry-eta-compact2}
|u_\epsilon(x_0)|\ge 1-L\eta^\theta.
\end{equation}
\end{theorem}

In order to prove Theorem \ref{bdry-eta-compact1}, we first establish the boundary monotonicity
inequalities \eqref{bdry_mono01.1} and \eqref{bdry_mono01.2} (see also \eqref{bdry_mono1.1} and \eqref{bdry_mono1.2})
for the quantity $\Phi_\epsilon(u_\epsilon, B_r(x_0)\cap\Omega)$ {with
$x_0\in\partial\Omega$. To do this, we employ a Pohozaev argument in $B_r(x_0)\cap\Omega$ by testing \eqref{GL-WA0} against
$X\cdot\nabla (u_\epsilon-\widetilde{g_\epsilon})$ with a certain vector field $X\in C^2(B_r(x_0),\mathbb R^n)$, satisfying $X(x)\in T_x(\partial\Omega)$ for $x\in\partial\Omega\cap B_r(x_0)$,}
and a suitable extension $\widetilde{g_\epsilon}$ of $g_\epsilon$. This involves some careful estimates on the error terms
$\int_{B_r(x_0)\cap\Omega}\frac{1-|u_\epsilon|^2}{\epsilon^2}u_\epsilon X\cdot\nabla\widetilde{g_\epsilon}$, to which we adopt an idea originated from \cite{CL} and
\cite{LR1}, and $\int_{B_r(x_0)\cap\Omega} (\nabla u_\epsilon\otimes\nabla\widetilde{g_\epsilon}:\nabla X
+\nabla u_\epsilon\otimes X:\nabla^2{\widetilde{g_\epsilon}})$.
With these boundary monotonicity inequalities, we then adopt the Hodge-decomposition techniques developed by \cite{BBO2} Appendix
into our setting of weak anchoring boundary conditions to ``clean out" any possible vortex near $x_0$, under the $\eta$-smallness condition
\eqref{eta-condition}.
To do so, we need to overcome several new difficulties arising from the equation
${\rm{div}}(\nabla u_\epsilon\times u_\epsilon)=0$ under the weak anchoring boundary condition \eqref{GL-WA0}$_2$. Among these difficulties,
we would like to mention that it is rather nontrivial to show $x\mapsto\int_{B_{r_0}(x_0^*)\cap\Omega}|x-y|^{2-n}\frac{(1-|u_\epsilon(y)|^2)^2}{\epsilon^2}\,dy$ is bounded in $B_{r_0}(x_0^*)$,
independent of $\epsilon$, where $x_0^*\in\partial\Omega$ satisfies $|x_0-x_0^*|={\rm{dist}}(x_0,\partial\Omega)$. We achieve this by employing both the interior
monotonicity inequality by \cite{BBO2} Lemma II.2 and the boundary monotonicity inequalities given by Theorem 2.3 and Theorem 2.4 {below. See Lemma 3.3 for the details.}

We would like to point out {that it is left open in} Theorem \ref{compactness1} whether the Hausdorff dimension of $\Sigma\cap\partial\Omega$
is $(n-3)$ and $u_*$ is H\"older continuous near $\partial\Omega\setminus\Sigma$
when $0<\alpha<1$. See Remark 6.3 below for more details.

{Combining the boundary monotonocity formula (\ref{bdry_mono1.1}) and the boundary $\eta$-compactness Theorem (\ref{bdry-eta-compact1}) with the interior monotonicity formula and the interior $\eta$-compactness
in \cite{BBO2} (Lemma II.2 and Theorem 2), we can show the following global property,} asserting
the uniform bound of potential energy over any approximate vortex set in $\overline\Omega$,
which plays a crucial role {in the proof of the $W^{1,p}$-estimate} for $\displaystyle 1\le p<\frac{n}{n-1}$.
More precisely, if we define the closed subset $S_\beta^\epsilon\subset\overline\Omega$, $0<\beta<\frac12$,  by
$$S^\epsilon_\beta:=\Big\{x\in\overline\Omega: \ \big|u_\epsilon(x)\big|\le 1-\beta\Big\},$$
then we have
\begin{theorem} \label{potential_bound0}
There exists $C_\beta>0$, depending on $\Omega$,
$\beta$, $K, \alpha$, $C_0$, and $M$,
such that for any $g_\epsilon\in C^2(\partial\Omega,\mathbb S^1)$ satisfying the condition $({\bf G})$,
if $u_\epsilon\in C^2(\overline\Omega,\R^2)$ is a solution {to
\eqref{GL-WA0} with $\lambda_\epsilon=K\epsilon^{-\alpha}$ for some
$K>0$ and $\alpha\in [0,1)$ that satisfies \eqref{global_bound}
for some $M>0$, then}
\beq\label{potential_bound00}
\int_{S_{\beta}^\epsilon}\frac{(1-|u_\epsilon|^2)^2}{\epsilon^2}\le C_\beta, \ \forall\epsilon\in (0,1].
\eeq
\end{theorem}

Utilizing Theorem \ref{potential_bound0}, we can extend the Hodge decomposition techniques
developed by \cite{BBO2} Appendix, with suitable nontrivial modifications due to the weak anchoring boundary condition,
to prove the global $W^{1,p}$-estimate that is a crucial part of Theorem \ref{compactness1}.
In contrast with \cite{BBO2}, we apply the Hodge decomposition directly to the weighted $1$-form
$\frac{1}{|du_\epsilon|^q}u_\epsilon\times du_\epsilon$ for some $q>0$,
{which seems to simplify} the whole argument substantially. More precisely, we have
\begin{theorem}\label{Lp-estimate0} For any $0<\epsilon\le 1$ and $g_\epsilon\in C^2(\partial\Omega,\mathbb S^1)$ satisfying the condition $({\bf G})$,
if $u_\epsilon\in C^2(\overline\Omega,\R^2)$ is a solution to
\eqref{GL-WA0}, with $\lambda_\epsilon=K\epsilon^{-\alpha}$ for some
$K>0$ and $\alpha\in [0,1)$, {that satisfies \eqref{global_bound}
for some $M>0$, then} for any $1\le p<\frac{n}{n-1}$,
there exists $C_p>0$, independent of $\epsilon$,  such that
\beq\label{Lp-estimate00}
\int_\Omega |\nabla u_\epsilon|^p\le C_p, \ \forall \epsilon\in (0,1].
\eeq
\end{theorem}

The paper is organized as follows. In section 2, we will derive two different versions of  boundary
monotonicity inequalities for \eqref{GL-WA0}. In section 3, we will prove the boundary $\eta$-compactness
Theorem \ref{bdry-eta-compact1}. In section 4, we will prove Theorem \ref{potential_bound0}.
In section 5, we will prove Theorem \ref{Lp-estimate0}.  Finally, in section 6, we will prove Theorem \ref{compactness1}.

\section{Boundary monotonicity inequality}
\setcounter{equation}{0}
\setcounter{theorem}{0}

In this section, we derive the boundary monotonicity inequalities for \eqref{GL-WA0} in dimensions $n\ge 3$.
For $x_0\in \R^n$ and $r>0$, let $B_r(x_0)$ denote the ball in $\R^n$ with center $x_0$ and radius
$r>0$.   For $x_0\in\overline\Omega$, set\footnote{If $r<{\rm{dist}}(x_0,\partial\Omega)$,  then
$B_r^+(x_0)=B_r(x_0)$, $S_r(x_0)=\partial B_r(x_0)$, and $\Gamma_r(x_0)=\emptyset$}
$$B_r^+(x_0)=B_r(x_0)\cap \Omega,\ S_r(x_0)=\partial B_r(x_0)\cap\Omega,
\ {\rm{and}}\ \ \Gamma_r(x_0)=\partial\Omega\cap B_r(x_0).$$

For $\epsilon>0$, $g_\epsilon\in C^\infty(\partial\Omega, \mathbb R^2)$, and $\lambda_\epsilon=K\epsilon^{-\alpha}$ with $K>0$ and $\alpha\in (0,1)$,
the Ginzburg-Landau equation with the weak anchoring boundary condition is given by
\begin{equation}\label{GZ1}\displaystyle
\begin{cases}
\Delta u_\epsilon+\frac{1}{\epsilon^2} (1-|u_\epsilon|^2)u_\epsilon =0 & \ {\rm{in}}\ \Omega,\\
\frac{\partial u_\epsilon}{\partial \nu}+\lambda_\epsilon(u_\epsilon-g_\epsilon)=0 & \ {\rm{on}}\ \partial\Omega.
\end{cases}
\end{equation}
Denote the Ginzburg-Landau energy density of $u_\epsilon$ by
$$e_\epsilon(u_\epsilon):=\frac{|\nabla u_\epsilon|^2}{2}+\frac{(1-|u_\epsilon|^2)^2}{4\epsilon^2}.$$
For $x_0\in\overline\Omega$ and $r>0$, the Ginzburg-Landau energy of $u_\epsilon$ on $B_r^+(x_0)$
is defined by
\begin{equation}\label{modified_energy1}
F_\epsilon^+(u_\epsilon; x_0, r):=\int_{B_r^+(x_0)}e_\epsilon(u_\epsilon)\,dx+\int_{\Gamma_r(x_0)}\frac{\lambda_\epsilon}2|u_\epsilon-g_\epsilon|^2\,dH^{n-1},
\end{equation}
so that
\begin{equation}
\label{derivative_modified_enerngy1}
\frac{d}{dr} F_\epsilon^+(u_\epsilon; x_0, r)=\int_{S_r(x_0)}e_\epsilon(u_\epsilon)\,dx+\int_{\partial\Gamma_r(x_0)}\frac{\lambda_\epsilon}2|u_\epsilon-g_\epsilon|^2\,dH^{n-2}
\end{equation}
for a.e. $r>0$.
Then we have
\begin{lemma}\label{stationarity1} Assume that the condition $({\bf G})$ holds.
Then there {exist $0<r_0=r_0(\Omega)<1$ and $C_1=C_1(\Omega)>0$}
such that if $u_\epsilon\in C^2\big(\overline\Omega, \mathbb R^2\big)$ is a solution of (\ref{GZ1}) then for any $x_0\in\partial\Omega$
and a.e. $0<r\le r_0$, {we have}
\begin{eqnarray}\label{stationarity2}
&&-C_1\Big[r^{n-1}+r F_\epsilon^+(u_\epsilon; x_0, r)\Big]
+\frac12 \Big[r\int_{S_r(x_0)}|\frac{\partial u_\epsilon}{\partial r}|^2+\int_{B_r^+(x_0)}\frac{(1-|u_\epsilon|^2)^2}{2\epsilon^2}\nonumber\\
&&\qquad+\int_{\Gamma_r(x_0)}\frac{\lambda_\epsilon}2 |u_\epsilon-g_\epsilon|^2\Big]
\nonumber\\
&&\le (2-n) F_\epsilon^+(u_\epsilon; x_0, r)+r\frac{d}{dr} F_\epsilon^+(u_\epsilon; x_0, r)\\
&&\le C_1\Big[r^{n-1}+r F_\epsilon^+(u_\epsilon; x_0, r)\Big]
+2\Big[r\int_{S_r(x_0)}|\frac{\partial u_\epsilon}{\partial r}|^2+\int_{B_r^+(x_0)}\frac{(1-|u_\epsilon|^2)^2}{2\epsilon^2}\nonumber\\
&&\qquad+\int_{\Gamma_r(x_0)}\frac{\lambda_\epsilon}2 |u_\epsilon-g_\epsilon|^2\Big].\nonumber
\end{eqnarray}
\end{lemma}
\pf Since $\partial\Omega$ is smooth,
there exist $r_0=r_0(\Omega)>0$ and $C_0=C_0(\Omega)>0$ such that for any $x_0\in\partial\Omega$,
there exists $X\in C^2(B_{r_0}(x_0),\mathbb R^n)$ {satisfying}
\begin{equation}\label{goodvf}
X\cdot\nu=0 \ {\rm{on}}\ T_{r_0}(x_0), \ |X(x)-(x-x_0)|\le C_0|x-x_0|^2, \ |DX(x)-\mathbb I_n|\le C_0|x-x_0|
\ {\rm{in}}\ B_{r_0}(x_0).
\end{equation}

It follows from the condition ({\bf G}) that there exists an extension $\widetilde{g_\epsilon}:B_{r_0}^+(x_0)\to\mathbb R^2$ of $g_\epsilon:T_{r_0}(x_0)\to\mathbb S^1$ such that
\begin{equation}\label{g-extension}
\|\widetilde{g_\epsilon}\|_{C^2(B_{r_0}^+(x_0))}\le 10C_0.
\end{equation}
To simplify the notation, we denote $u=u_\epsilon, \lambda=\lambda_\epsilon,$ $g=g_\epsilon$, and $\widetilde{g}=\widetilde{g_\epsilon}$.
For $0<r\le r_0$, multiplying (\ref{GZ1})$_1$ by $X\cdot \nabla (u-\widetilde{g})$, integrating the resulting equation in $B_r^+(x_0)$, and applying integration by parts, we obtain
\begin{eqnarray} \label{stationarity3}
0&=&\int_{B_r^+(x_0)}\nabla\cdot\langle X\cdot \nabla (u-\widetilde{g}), \nabla u\rangle
-\nabla u\otimes\nabla (u-\widetilde{g}):\nabla X\nonumber\\
&&+\int_{B_r^+(x_0)}X\otimes\nabla u:\nabla^2 \widetilde{g}-\int_{B_r^+(x_0)}X\cdot\nabla(e_\epsilon(u))
-\int_{B_r^+(x_0)}\frac{1-|u|^2}{\epsilon^2} u X\cdot\nabla\widetilde{g}\nonumber\\
&=&\int_{S_r(x_0)}\big(\langle\frac{\partial u}{\partial r}, X\cdot\nabla u\rangle-e_\epsilon(u)\frac{X\cdot (x-x_0)}{r}\big)
+\int_{B_r^+(x_0)}\big[e_\epsilon(u){\rm{div}}X-\nabla u\otimes\nabla u:\nabla X\big]\nonumber\\
&&+\int_{\Gamma_r(x_0)} \langle \frac{\partial u}{\partial\nu}, X\cdot\nabla (u-\widetilde{g})\rangle
-\int_{B_r^+(x_0)}\frac{1-|u|^2}{\epsilon^2} u X\cdot\nabla\widetilde{g}\nonumber\\
&&+\int_{B_r^+(x_0)}\big[\nabla u\otimes\nabla \widetilde{g}:\nabla X+X\otimes\nabla u:\nabla^2 \widetilde{g}\big]
-\int_{S_r(x_0)}\langle\frac{\partial u}{\partial r}, X\cdot\nabla \widetilde{g}\rangle.
\end{eqnarray}
From \eqref{goodvf}, we have that $X(x)\in T_x\Gamma_r(x_0)$ for $x\in\Gamma_r(x_0)$.
Hence, using (\ref{GZ1})$_2$, we can calculate
\begin{eqnarray}\label{estimate0.1}
&&\int_{\Gamma_r(x_0)}\langle \frac{\partial u}{\partial\nu}, X\cdot\nabla (u-\widetilde{g})\rangle=-\lambda
\int_{\Gamma_r(x_0)}\langle X\cdot\nabla (u-\widetilde{g}), u-\widetilde{g}\rangle\nonumber\\
&&=-\lambda
\int_{\Gamma_r(x_0)}{\rm{div}}_{\Gamma_{r}(x_0)}\big(\frac{|u-\widetilde{g}|^2}2 X\big)
+\lambda\int_{\Gamma_r(x_0)}\frac{|u-\widetilde{g}|^2}2
{\rm{div}}_{\Gamma_{r}(x_0)}(X)
\nonumber\\
&&=-
\int_{\partial\Gamma_r(x_0)}\frac{\lambda}2 |u-{g}|^2X\cdot \nu_{\partial\Gamma_r(x_0)}(x)
+\int_{\Gamma_r(x_0)}\frac{\lambda}2 |u-{g}|^2{\rm{div}}_{\Gamma_{r}(x_0)} (X).
\end{eqnarray}

Note that it follows from \eqref{crude_bound} below that $|u|\le 1$ in $B_r^+(x_0)$.
Now we proceed as follows.

\smallskip
The strategy to estimate  $\int_{B_r^+(x_0)}\frac{1-|u|^2}{\epsilon^2} u X\cdot\nabla\widetilde{g}$ is similar to that of Lin-Rivi\`ere \cite{LR1} and Chen-Lin \cite{CL}. More precisely,
let $\phi\in C^\infty([0,1],[0,1])$ be such that $\phi(0)=0$, $\phi(t)=1$ for $\epsilon^2\le t\le 1$,
and $\phi'(t)\ge 0$ for $t\in [0,1]$. Multiplying equation \eqref{GZ1}$_1$ by $\phi(1-|u|^2)u$ and integrating over
$B_r^+(x_0)$, we obtain
\begin{eqnarray}\label{estimate0}
&&\int_{B_r^+(x_0)}\frac{(1-|u|^2)|u|^2}{\epsilon^2}\phi(1-|u|^2)\nonumber\\
&&=\int_{B_r^+(x_0)}\nabla u\cdot\nabla(\phi(1-|u|^2)u)-\int_{\partial B_r^+(x_0)}\frac{\partial u}{\partial\nu} \phi(1-|u|^2)u\nonumber\\
&&=\int_{B_r^+(x_0)}|\nabla u|^2\phi(1-|u|^2)-2\int_{B_r^+(x_0)}\phi'(1-|u|^2)|\langle u,\nabla u\rangle|^2\nonumber\\
&&\quad-\int_{S_r(x_0)}\frac{\partial u}{\partial r} \phi(1-|u|^2)u+\lambda\int_{T_r(x_0)}\langle u, u-g\rangle\phi(1-|u|^2)
\nonumber\\
&&\le \int_{B_r^+(x_0)}|\nabla u|^2+\int_{S_r(x_0)}|\frac{\partial u}{\partial r}|
+\lambda\int_{T_r(x_0)}\langle u, u-g\rangle\phi(1-|u|^2).
\end{eqnarray}
We claim that
\begin{equation}\label{estimate1}
\lambda\int_{T_r(x_0)}\langle u, u-g\rangle \phi(1-|u|^2)\le C\int_{T_r(x_0)}\lambda |u-g|^2.
\end{equation}
In fact, since $|g(x)|=1$ and $|u(x)|\le 1$ for any $x\in T_r(x_0)$, we have
$$|g|^2-u\cdot g\ge |g|^2-|u||g|=1-|u|\ge 0\ {\rm{on}}\ T_r(x_0),$$
and hence
\begin{eqnarray*}&&\langle u, u-g\rangle=|u|^2-u\cdot g\le |u|^2+|g|^2-2u\cdot g
=|u-g|^2 \ {\rm{on}}\ T_r(x_0).
\end{eqnarray*}
Thus we obtain
$$\lambda\int_{T_r(x_0)}\langle u, u-g\rangle \phi(1-|u|^2)\le \int_{T_r(x_0)}\lambda |u-g|^2.$$
Substituting \eqref{estimate1} into \eqref{estimate0},  we obtain
\begin{equation}\label{estimate1.1}
\int_{B_r^+(x_0)}\frac{(1-|u|^2)|u|^2}{\epsilon^2}\phi(1-|u|^2)
\le C\big(F^+_\epsilon(u; x_0, r)+\int_{S_r(x_0)}|\frac{\partial u}{\partial r}|\big).
\end{equation}
Applying \eqref{estimate1.1}, \eqref{goodvf} and \eqref{g-extension}, it is not hard to see that
\begin{eqnarray}\label{estimate1.2}
&&\big|\int_{B_r^+(x_0)}\frac{1-|u|^2}{\epsilon^2} u X\cdot\nabla\widetilde{g}\big|\le Cr\int_{B_r^+(x_0)}\frac{1-|u|^2}{\epsilon^2}\nonumber\\
&&\le  Cr\big(\int_{B_r^+(x_0)}\frac{(1-|u|^2)^2}{\epsilon^2}+\int_{B_r^+(x_0)}\frac{(1-|u|^2)}{\epsilon^2}|u|^2\big)\nonumber\\
&&\le Cr\big(\int_{B_r^+(x_0)}\frac{(1-|u|^2)^2}{\epsilon^2}+r^n
+\int_{B_r^+(x_0)}\frac{(1-|u|^2)}{\epsilon^2}|u|^2\phi(1-|u|^2)\big)\nonumber\\
&&\le Cr\big(r^{n}+F^+_\epsilon(u; x_0, r)+
\int_{S_r(x_0)}|\frac{\partial u}{\partial r}|\big)\nonumber\\
&&\le Cr^n+CrF^+_\epsilon(u; x_0, r)+\frac{1}{8}
r\int_{S_r(x_0)}|\frac{\partial u}{\partial r}|^2.
\end{eqnarray}
With the help of \eqref{g-extension}, we can estimate the last two terms of the right hand side of \eqref{stationarity3} as follows.
\begin{eqnarray}\label{estimate1.3}
\Big|\int_{B_r^+(x_0)}\big[\nabla u\otimes\nabla \widetilde{g}:\nabla X+X\otimes\nabla u:\nabla^2 \widetilde{g}\big]\Big|
&\le& C\int_{B_r^+(x_0)}|\nabla u|\nonumber\\
&\le& C\big(r^{n-1}+r\int_{B_r^+(x_0)}|\nabla u|^2\big),
\end{eqnarray}
and
\begin{eqnarray}\label{estimate1.4}
\big|-\int_{S_r(x_0)}\langle\frac{\partial u}{\partial r}, X\cdot\nabla \widetilde{g}\rangle\big|
&\le& Cr\int_{S_r(x_0)}|\frac{\partial u}{\partial r}|\nonumber\\
&\le& \frac1{8} r\int_{S_r(x_0)}|\frac{\partial u}{\partial r}|^2+Cr^{n}.
\end{eqnarray}
From \eqref{goodvf}, we have that
$$
\begin{cases} |X(x)|\le Cr, \ |{\rm{div}}X(x)-n|\le Cr, \ \forall x\in B_r^+(x_0),\\
|X\cdot\nu_{\partial\Gamma_r(x_0)}-r|\le Cr^2 \ {\rm{on}}\ \partial\Gamma_r(x_0),\\
|{\rm{div}}_{\Gamma_r(x_0)}(X)-(n-1)|\le Cr \ {\rm{on}}\ \Gamma_r(x_0).
\end{cases}
$$
Hence we can estimate
\begin{eqnarray}\label{bdry-term1}
&&\big|\int_{B_r^+(x_0)}\big[e_\epsilon(u){\rm{div}}X-\nabla u\otimes\nabla u:\nabla X\big]
-n\int_{B_r^+(x_0)}e_\epsilon(u)+\int_{B_r^+(x_0)}|\nabla u|^2\big|\nonumber\\
&&\le Cr\int_{B_r^+(x_0)}e_\epsilon(u),
\end{eqnarray}
\begin{eqnarray}\label{bdry-term2}
&&\big|\int_{S_r(x_0)}\langle \frac{\partial u}{\partial\nu}, X\cdot\nabla u\rangle-\int_{S_r(x_0)}e_\epsilon(u)\frac{X\cdot (x-x_0)}{r}
-r\int_{S_r(x_0)}|\frac{\partial u}{\partial r}|^2+r\int_{S_r(x_0)}e_\epsilon(u)\big|\nonumber\\
&&\le Cr^2\int_{S_r(x_0)}e_\epsilon(u),
\end{eqnarray}
and
\begin{eqnarray}\label{bdry-term3}
&&\Big|\int_{\Gamma_r(x_0)}\langle \frac{\partial u}{\partial\nu}, X\cdot\nabla (u-\widetilde{g})\rangle
+r\int_{\partial \Gamma_r(x_0)}\frac{\lambda}2|u-g|^2-(n-1)\int_{\Gamma_r(x_0)}\frac{\lambda}2 |u-g|^2\Big|\nonumber\\
&&\le Cr\int_{\Gamma_r(x_0)}\frac{\lambda}2 |u-g|^2+Cr^2\int_{\partial \Gamma_r(x_0)}\frac{\lambda}2|u-g|^2.
\end{eqnarray}
Putting \eqref{bdry-term1}, \eqref{bdry-term2}, \eqref{bdry-term3} into
\eqref{stationarity3}, we obtain that
\begin{eqnarray}\label{total_estimate1}
&&\Big|(2-n)F^+_\epsilon(u;x_0,r)+r\frac{d}{dr}F^+_\epsilon(u;x_0,r)-\big(\int_{B_r^+(x_0)}\frac{(1-|u|^2)^2}{2\epsilon^2}+\int_{\Gamma_r(x_0)}\frac{\lambda}2|u-g|^2\big)\nonumber\\
&&\ -r\int_{S_r(x_0)}|\frac{\partial u}{\partial r}|^2\Big|\nonumber\\
&&\le C\big[rF^+_\epsilon(u;x_0,r)+r^2\frac{d}{dr}F^+_\epsilon(u;x_0,r)+r^{n-1}\big]
+\frac18 r\int_{S_r(x_0)}|\frac{\partial u}{\partial r}|^2.
\end{eqnarray}
It follows from \eqref{total_estimate1} that
\begin{eqnarray*}
&&\frac{2-n}{1+Cr}F^+_\epsilon(u;x_0,r)+r\frac{d}{dr}F^+_\epsilon(u;x_0,r)\\
&&\ge \frac12\big(\int_{B_r^+(x_0)}\frac{(1-|u|^2)^2}{2\epsilon^2}+\int_{\Gamma_r(x_0)}\frac{\lambda}2|u-g|^2
+\int_{S_r(x_0)}|\frac{\partial u}{\partial r}|^2\big)-C\big(rF^+_\epsilon(u;x_0,r)+r^{n-1}\big),
\end{eqnarray*}
which implies the first inequality of \eqref{stationarity2}, since there exists $C_1>C$ such that $\frac{2-n}{1+Cr}\le 2-n+C_1r$.
Similarly, the second inequality of \eqref{stationarity2} can be obtained. \qed

\medskip
Now we deduce a few consequences of  (\ref{stationarity2}), that will be used in later sections.
The first one provides control of the tangential energy of $u_\epsilon$ on $S_r(x_0)$
in terms of both the radial energy of $u_\epsilon$ on $S_r(x_0)$ and the Ginzburg-Landau energy of $u_\epsilon$ in $B_r^+(x_0)$.

\begin{proposition}\label{tangential_radial1} Under the same assumptions as in Lemma \ref{stationarity1}, there exists $C_1>0$ such that
for any $x_0\in\partial\Omega$ and $0<r\le r_0$, we have
\beq\label{tangential_radial2}
r\int_{S_r(x_0)}\big|\nabla_T u_\epsilon\big|^2\le C\Big[r\int_{S_r(x_0)}\big|\frac{\partial u_\epsilon}{\partial |x-x_0|}\big|^2+\int_{B_r^+(x_0)}e_\epsilon(u_\epsilon)\Big]+C_1r^{n-1}.
\eeq
Here $\nabla_Tu_\epsilon$ denotes the tangential gradient of $u_\epsilon$ on $S_r(x_0)$.
\end{proposition}
\pf Since $\displaystyle |\nabla u_\epsilon|^2=|\nabla_T u_\epsilon|^2+\big|\frac{\partial u_\epsilon}{\partial |x-x_0|}\big|^2$ on $S_r(x_0)$, we have
$$\int_{S_r(x_0)}e_\epsilon(u_\epsilon)\ge \frac12\int_{S_r(x_0)}|\nabla u_\epsilon|^2=\frac12\int_{S_r(x_0)}|\nabla_T u_\epsilon|^2+
\frac12\int_{S_r(x_0)}\big|\frac{\partial u_\epsilon}{\partial |x-x_0|}\big|^2.$$
Substituting this into the second inequality of (\ref{stationarity2}), we can easily obtain (\ref{tangential_radial2}).
\qed

\medskip
Integrating (\ref{stationarity2}), we obtain two slightly different forms of boundary monotonicity inequalities
for $u$. The first one involves the renormalized Ginzburg-Landau energy of $u_\epsilon$ on $B_r^+(x_0)$ defined by
$$\Phi_\epsilon^+(u_\epsilon; x_0, r):=r^{2-n}F_\epsilon^+(u_\epsilon; x_0, r).$$

\begin{theorem}\label{bdry_mono00}
Assume the condition $({\bf G})$ holds.
Then there exist $0<r_0=r_0(\Omega)<1$ and $C_2=C_2(\Omega)>0$
such that if $u_\epsilon\in C^2\big(\overline\Omega, \mathbb R^2\big)$ is a solution of (\ref{GZ1}), then for any $x_0\in\partial\Omega$
and $0<r\le R\le r_0$, we have
\begin{eqnarray}\label{bdry_mono01.1}
&&e^{-C_2R}\Phi_\epsilon^+(u_\epsilon; x_0, R)
-e^{-C_2 r}\Phi_\epsilon^+(u_\epsilon; x_0, r)\nonumber\\
&&\le 2\Big[\int_{r}^R\tau^{1-n}\big(\int_{B_\tau^+(x_0)}\frac{(1-|u_\epsilon|^2)^2}{2\epsilon^2}+\int_{\Gamma_\tau(x_0)}\frac{\lambda_\epsilon}2|u_\epsilon-g_\epsilon|^2\big)\nonumber\\
&&\ \ +\int_{B_R^+(x_0)\setminus B_r^+(x_0)} |x-x_0|^{2-n}\big|\frac{\partial u_\epsilon}{\partial |x-x_0|}\big|^2\Big]+C_2(R-r),
\end{eqnarray}
and
\begin{eqnarray}\label{bdry_mono01.2}
&&\big(e^{C_2 R}\Phi_\epsilon^+(u_\epsilon; x_0, R)+C_2R\big)
-\big(e^{C_2 r}\Phi_\epsilon^+(u_\epsilon; x_0, r)+C_2r\big)\nonumber\\
&&\ge \frac12\Big[\int_{r}^R\tau^{1-n}\big(\int_{B_\tau^+(x_0)}\frac{(1-|u_\epsilon|^2)^2}{2\epsilon^2}
+\int_{\Gamma_\tau(x_0)}\frac{\lambda_\epsilon}2|u_\epsilon-g_\epsilon|^2\big)
\nonumber\\
&&\qquad +\int_{B_R^+(x_0)\setminus B_r^+(x_0)} |x-x_0|^{2-n}\big|\frac{\partial u_\epsilon}{\partial |x-x_0|}\big|^2\Big].
\end{eqnarray}
\end{theorem}
\pf It follows directly from the first inequality of (\ref{stationarity2}) that
\begin{eqnarray*}
&&\frac{d}{dr}\big(r^{2-n}F_\epsilon^+(u_\epsilon; x_0, r)\big)+C_1r^{2-n}F_\epsilon^+(u_\epsilon; x_0, r)+C_1\\
&&\ge \frac12r^{1-n}\Big(r\int_{S_r(x_0)}\big|\frac{\partial u_\epsilon}{\partial |x-x_0|}\big|^2+\int_{B_r^+(x_0)}\frac{(1-|u_\epsilon|^2)^2}{2\epsilon^2}
+\int_{\Gamma_r(x_0)}\frac{\lambda_\epsilon}{2}|u_\epsilon-g_\epsilon|^2\Big).
\end{eqnarray*}
Integrating this inequality over $[r, R]$ yields (\ref{bdry_mono01.2}). Similarly, it follows from the second inequality of \eqref{stationarity2}
that
\begin{eqnarray*}
&&\frac{d}{dr}\big(r^{2-n}F_\epsilon^+(u_\epsilon; x_0, r)\big)-C_1r^{2-n}F_\epsilon^+(u_\epsilon; x_0, r)\\
&&\le 2r^{1-n}\Big(r\int_{S_r(x_0)}\big|\frac{\partial u_\epsilon}{\partial |x-x_0|}\big|^2+\int_{B_r^+(x_0)}\frac{(1-|u_\epsilon|^2)^2}{2\epsilon^2}
+\int_{\Gamma_r(x_0)}\frac{\lambda_\epsilon}{2}|u_\epsilon-g_\epsilon|^2\Big)+C_1
\end{eqnarray*}
which, after integrating over $[r,R]$, implies (\ref{bdry_mono01.1}).
\qed

\medskip
To state the second form of the boundary monotonicity inequality, we {define another form of a modified Ginzburg-Landau energy density of}
$u_\epsilon$ by
$$\widehat{e}_\epsilon(u_\epsilon):=\frac{|\nabla u_\epsilon|^2}2+\frac{n(1-|u_\epsilon|^2)^2}{2(n-2)\epsilon^2},$$
so that the corresponding modified Ginzburg-Landau energy of $u_\epsilon$ in $B_{r}^+(x_0)$, {for $x_0\in\overline\Omega$ and $r>0$, is}

$$
\widehat{F}_\epsilon^+(u_\epsilon; x_0, r)=\int_{B_r^+(x_0)}\widehat{e}_\epsilon(u_\epsilon)
+\int_{\Gamma_r(x_0)} \frac{(n-1)\lambda_\epsilon |u_\epsilon-g_\epsilon|^2}{2(n-2)}\,dH^{n-1},
$$
and the corresponding {red}{renormalized energy is}
$$
\widehat{\Phi}_\epsilon^+(u_\epsilon; x_0,r)=r^{2-n}\widehat{F}_\epsilon^+(u_\epsilon; x_0, r).
$$
Since
\begin{equation}\label{comparable_density}
e_\epsilon(u_\epsilon)\le \widehat{e}_\epsilon(u_\epsilon)\le \frac{2n}{n-2} e_\epsilon(u_\epsilon),
\end{equation}
it is easy to see that
\begin{equation}\label{comparble_enegry}
\Phi^+_\epsilon(u_\epsilon; x_0, r)\le \widehat{\Phi}_\epsilon^+(u_\epsilon; x_0, r)
\le \frac{2n}{n-2}\Phi^+_\epsilon(u_\epsilon; x_0, r).
\end{equation}
Then we have
\begin{theorem}\label{bdry_mono0} Assume that the condition $({\bf G})$ holds.
Then there exist $0<r_0=r_0(\Omega)<1$ and $C_3=C_3(\Omega)>0$
such that if $u_\epsilon\in C^2\big(\overline\Omega, \mathbb R^2\big)$ is a solution of (\ref{GZ1}), then for any $x_0\in\partial\Omega$
and $0<r\le R\le r_0$, we have
\begin{eqnarray}\label{bdry_mono1.1}
&&e^{-C_3 R}\widehat{\Phi}_\epsilon^+(u_\epsilon; x_0, R)-e^{-C_3 r} \widehat{\Phi}_\epsilon^+(u_\epsilon; x_0, r)\nonumber\\
&&\le 2\Big[\int_{B_R^+(x_0)\setminus B_r^+(x_0)}|x-x_0|^{2-n}\big(|\frac{\partial u_\epsilon}{\partial |x-x_0|}|^2+\frac{(1-|u_\epsilon|^2)^2}{2(n-2)\epsilon^2}\big)\nonumber\\
&&\ \ \ \ \ \ +\int_{\Gamma_R(x_0)\setminus\Gamma_r(x_0)} |x-x_0|^{2-n}\frac{\lambda_\epsilon |u_\epsilon-g_\epsilon|^2}{2(n-2)}\Big]+C_3(R-r),
\end{eqnarray}
and
\begin{eqnarray}\label{bdry_mono1.2}
&&e^{C_3 R}\widehat{\Phi}_\epsilon^+(u_\epsilon; x_0, R)- e^{C_3r}\widehat{\Phi}_\epsilon^+(u_\epsilon; x_0, r)+C_3(R-r)\nonumber\\
&&\ge \frac12\Big[\int_{B_R^+(x_0)\setminus B_r^+(x_0)}|x-x_0|^{2-n}\Big(\big|\frac{\partial u_\epsilon}{\partial |x-x_0|}\big|^2+\frac{(1-|u_\epsilon|^2)^2}{2(n-2)\epsilon^2}\Big)\nonumber\\
&&+\int_{\Gamma_R(x_0)\setminus\Gamma_r(x_0)} |x-x_0|^{2-n}\frac{\lambda_\epsilon |u_\epsilon-g_\epsilon|^2}{2(n-2)}\Big].
\end{eqnarray}
\end{theorem}
\pf It is easy to see that the first inequality (\ref{stationarity2}) implies that
\begin{eqnarray*}
&&r\big[\int_{S_r(x_0)}\widehat{e}_\epsilon(u_\epsilon)+\int_{\partial\Gamma_r(x_0)}\frac{(n-1)\lambda_\epsilon |u_\epsilon-g_\epsilon|^2}{2(n-2)}\big]-C_1r^{n-1}\\
&&\le (n-2+C_1 r)
\big[\int_{B_r^+(x_0)}\widehat{e}_\epsilon(u_\epsilon)+\int_{\Gamma_r(x_0)}\frac{(n-1)\lambda_\epsilon |u_\epsilon-g_\epsilon|^2}{2(n-2)}\big]\\
&&+2\Big[r\int_{S_r(x_0)}\big[\big|\frac{\partial u_\epsilon}{\partial |x-x_0|}\big|^2+\frac{(1-|u_\epsilon|^2)^2}{2(n-2)\epsilon^2}\big]
+r\int_{\partial\Gamma_r(x_0)}\frac{\lambda_\epsilon |u_\epsilon-g_\epsilon|^2}{2(n-2)}\Big].
\end{eqnarray*}
Thus we have that
\begin{eqnarray}
\label{bdry_mono3}
&&\frac{d}{dr}\Big[e^{-C_3 r}r^{2-n}\big(\int_{B_r^+(x_0)}\widehat{e}_\epsilon(u_\epsilon)+\int_{\Gamma_r(x_0)}\frac{(n-1)\lambda_\epsilon |u_\epsilon-g_\epsilon|^2}{2(n-2)}\big)\Big]
-C_3\nonumber\\
&&\le 2r^{2-n}\Big[\int_{S_r(x_0)}\big(|\frac{\partial u_\epsilon}{\partial |x-x_0|}|^2+\frac{(1-|u_\epsilon|^2)^2}{2(n-2)\epsilon^2}\big)
+\int_{\partial\Gamma_r(x_0)} \frac{\lambda_\epsilon |u_\epsilon-g_\epsilon|^2}{2(n-2)}\Big].
\end{eqnarray}
Integrating (\ref{bdry_mono3}) from $r$ to $R$ yields (\ref{bdry_mono1.1}). Similarly, (\ref{bdry_mono1.2})
can be derived by integrating the second inequality of \eqref{stationarity2}.
\qed

\medskip
We can draw an immediate consequence of Theorem \ref{bdry_mono0} as follows.
\begin{corollary} Assume that the condition $({\bf G})$ holds.
Then there exist $0<r_0=r_0(\Omega)<1$ and $C_4=C_4(\Omega)>0$
such that if $u_\epsilon\in C^2\big(\overline\Omega, \mathbb R^2\big)$ is a solution of (\ref{GZ1}), then for any $x_0\in\partial\Omega$
and $0<R\le r_0$, we have
\begin{eqnarray}\label{bdry_mono5}
&&\int_{B_R^+(x_0)}|x-x_0|^{2-n}\Big(\big|\frac{\partial u_\epsilon}{\partial |x-x_0|}\big|^2+\frac{(1-|u_\epsilon|^2)^2}{2(n-2)\epsilon^2}\Big)
+\int_{\Gamma_R(x_0)}|x-x_0|^{2-n}\frac{\lambda_\epsilon|u_\epsilon-g_\epsilon|^2}{2(n-2)}\nonumber\\
&&\le e^{ C_4R}\widehat{\Phi}^+_\epsilon(u_\epsilon; x_0, R)+C_4 R.
\end{eqnarray}
\end{corollary}
\pf It is easy to see that (\ref{bdry_mono5}) follows from (\ref{bdry_mono1.2}) by sending $r\rightarrow 0$.
\qed

\medskip
Similar to \cite{BBO2} and \cite{LR1}, in order to prove Theorem \ref{potential_bound0} we also need the following bound on the renormalized Ginzburg-Landau energy, namely:

\begin{theorem}\label{RGLB} Under the same assumptions as Theorem \ref{compactness1}, {for any $x_0\in\overline\Omega$, we have}
\begin{equation}\label{RGLB1}
\Phi^+_\epsilon(u_\epsilon; x_0, r)=r^{2-n}\big(\int_{B_r^+(x_0)}e_\epsilon(u_\epsilon)+\int_{\Gamma_r(x_0)}\frac{\lambda_\epsilon}2|u_\epsilon-g_\epsilon|^2\big)\le C|\log\epsilon|, \ \forall r>0,
\end{equation}
where $C>0$ is a constant that is independent of $x_0, r$, and $\epsilon$.
\end{theorem}
\pf This follows directly from \eqref{global_bound}, \eqref{bdry_mono01.1}, and the interior energy monotonicity inequality. (See \cite{BBO2}) and \cite{LR1} for $u_\epsilon$.)
\qed

\section{The boundary $\eta$-compactness property}
\setcounter{equation}{0}
\setcounter{theorem}{0}

In this section, we will derive the boundary $\eta$-compactness property for solutions to \eqref{GL-WA0}
under the weak anchoring boundary condition.

From now on, we assume that $\lambda=\lambda_\epsilon\equiv
K\epsilon^{-\alpha}$ for some $K>0$ and $\alpha\in [0,1)$.
First we need to establish the following crude estimate:
\begin{lemma} For any $g_\epsilon\in C^2(\partial\Omega, \mathbb S^1)$, there exists
$C>0$ depending on $K,\alpha$, $\Omega$, and $\|g_\epsilon\|_{C^1(\partial\Omega)}$
such that if $u_\epsilon\in C^2(\overline\Omega,\mathbb R^2)$
is a solution of
\begin{equation}\label{GL_WA1}
\begin{cases} \Delta u_\epsilon+\frac{1}{\epsilon^2}(1-|u_\epsilon|^2)u_\epsilon =0
 & {\rm{in}}\ \Omega,\\
\frac{\partial u_\epsilon}{\partial\nu}+\lambda_\epsilon (u_\epsilon-g_\epsilon)=0 & {\rm{on}}\ \partial\Omega,
\end{cases}
\end{equation}
then
\begin{equation}\label{crude_bound}
|u_\epsilon(x)|\le 1, \ |\nabla u_\epsilon(x)|\le C\epsilon^{-1}, \ \forall\ x\in\overline\Omega.
\end{equation}
\end{lemma}
\pf Set $W_\epsilon=|u_\epsilon|^2-1$ and $W_\epsilon^+=\max\big\{W_\epsilon, 0\big\}$.
Then it is easy to check that
$$\Delta W_\epsilon=\frac{2}{\epsilon^2}(W_\epsilon+1)W_\epsilon+2|\nabla u_\epsilon|^2
\ge \frac{2}{\epsilon^2}(W_\epsilon+1)W_\epsilon, \ {\rm{in}}\ \ \Omega.$$
Multiplying this equation by $W_\epsilon^+$ and integrating over $\Omega$, we obtain
$$\frac{2}{\epsilon^2}\int_\Omega |u_\epsilon|^2 (W_\epsilon^+)^2
\le \int_{\partial\Omega} W_\epsilon^+\frac{\partial W_\epsilon}{\partial\nu}
-\int_\Omega |\nabla W_\epsilon^+|^2.$$
From \eqref{GL_WA1}$_2$, we have that for any $x\in\partial\Omega$,
\begin{eqnarray*}
W_\epsilon^+(x)\frac{\partial W_\epsilon}{\partial\nu}(x)
&=&-2\lambda_\epsilon W_\epsilon^+(x)u_\epsilon(x)\cdot (u_\epsilon(x)-g_\epsilon(x))\\
&=&-2\lambda_\epsilon W_\epsilon^+(x)(|u_\epsilon(x)|^2-u_\epsilon(x)\cdot g_\epsilon(x))\\
&\le& -2\lambda_\epsilon W_\epsilon^+(x)(|u_\epsilon(x)|^2-|u_\epsilon(x)|)\le 0.
\end{eqnarray*}
Thus we obtain
$$\frac{2}{\epsilon^2}\int_\Omega |u_\epsilon|^2 (W_\epsilon^+)^2\le 0,$$
which implies that $W_\epsilon^+\equiv 0$ and hence $|u_\epsilon(x)|\le 1$ on $\overline\Omega$.

The gradient estimate of $u_\epsilon$ can be proved by a contradiction argument. Suppose it were false. Then
there exist $\epsilon_k\rightarrow 0$ and a sequence $u_k(=u_{\epsilon_k})\in C^2(\overline\Omega,\R^2)$ solving (\ref{GL_WA1}) such that
$m_k=\big\|\nabla u_k\big\|_{L^\infty(\overline\Omega)}$ satisfies $m_k\epsilon_k\rightarrow\infty$.
Let $x_k\in\overline\Omega$ be such that $|\nabla u_k(x_k)|=m_k$. Define the blow-up sequence
{$$v_k(x)=u_k(x_k+\frac{x}{m_k}\big), \ \text{for } x\in\Omega_k=m_k(\Omega - \{x_k\}).$$}
Then
\begin{equation}\label{GL_WA2}
\begin{cases}
\Delta v_k+\frac{1}{(m_k\epsilon_k)^2}(1-|v_k|^2)v_k=0 & {\rm{in}}\ \Omega_k,\\
\frac{\partial v_k}{\partial \nu}+\frac{\lambda_k}{m_k}(v_k-g_k)=0 & {\rm{on}}\ \partial\Omega_k,
\end{cases}
\end{equation}
where $g_k(x)\equiv g_{\epsilon_k}\big(x_k+\frac{x}{m_k}\big)$ for $x\in\partial\Omega_k$ {and $\lambda_k=\lambda_{\epsilon_k}$.}  Moreover,
\begin{equation}\label{total_bound}
|\nabla v_k(0)|=1=\big\|\nabla v_k\|_{L^\infty(\overline\Omega_k)}, \ \big\|v_k\big\|_{L^\infty(\overline\Omega_k)}\le 1.
\end{equation}
Now we divide the proof into two cases.\\
(i) $m_k{\rm{dist}}(x_k,\partial\Omega)\rightarrow\infty$ and $\Omega_k\rightarrow \R^n$. It follows from (\ref{total_bound})
that
$$\max_k\|v_k\|_{C^l(B_R)}\le C(R,l)<+\infty$$
for any $l\ge 1$ and $0<R<+\infty$. Therefore we can assume
that $v_k\rightarrow v_\infty$ in $C^2_{\rm{loc}}(\R^n)$. It is clear from (\ref{GL_WA2}) that $v_\infty$ satisfies
$$\Delta v_\infty=0 \ {\rm{in}}\ \R^n, $$
$$ |\nabla v_\infty(0)|=1 \ {\rm{and}}\ \|v_\infty\|_{L^\infty(\R^n)}\le 1.$$
This is impossible.\\
(ii) $m_k{\rm{dist}}(x_k,\partial\Omega)\rightarrow a$ for some $0\le a<+\infty$,
and $\Omega_k\rightarrow \R^n_{-a}:=\big\{x=(x', x_n)\in\R^n: x_n\ge -a\big\}$.
Since $\frac{\lambda_k}{m_k}=\frac{K\epsilon_k^{1-\alpha}}{m_k\epsilon_k}\rightarrow 0$,
it follows from (\ref{GL_WA2}) and
(\ref{total_bound}) that
$$\max_k\big\|v_k\big\|_{C^l(B_R\cap \Omega_k)}\le C(R, l), \ \forall\ l\ge 1, \  R>0.$$
Thus we may assume that $v_k\rightarrow v_\infty$ in $C^2_{\rm{loc}}(\R^n_{-a})$ so that
$v_\infty$ solves
$$\Delta v_\infty=0 \ {\rm{in}}\ \R^n_{-a}; \ \frac{\partial v_\infty}{\partial \nu}=0 \ {\rm{on}}\ \partial \R^n_{-a},$$
and
$$|\nabla v_\infty(0)|=1\ {\rm{and}} \ \big\|v_k\big\|_{L^\infty(\R^n_{-a})}\le 1.$$
This is again impossible. The proof is now complete.
\qed

\medskip

We are ready to prove the following global $\eta$-compactness property for \eqref{GL_WA1}.

\begin{theorem} \label{eta-comp} There exist $\epsilon_0>0$, depending on $\eta$,
$\Omega, C_0, K,\alpha$, and $r_0=r_0(\Omega)>0$ such that
for any $\eta>0$, $K>0$, $\alpha\in [0,1)$, if $g_\epsilon\in C^2(\partial\Omega,\mathbb S^1)$ satisfies the condition
$({\bf G})$ and $u_\epsilon\in C^2(\overline{\Omega},\mathbb R^2)$ is a solution of (\ref{GL_WA1}), with $\lambda_\epsilon=K\epsilon^{-\alpha}$, then
the following is true: if for a fixed $x_0\in\overline\Omega$
\begin{equation}\label{eta-com1}
\Phi^+_\epsilon(u_\epsilon; x_0, r)
\le\eta |\log \epsilon|
\end{equation}
holds for $0<\epsilon\le\epsilon_0$ and $\epsilon^{\alpha}\le r\le r_0$,
then there exist two constants $L>0$ and $\theta\in (0,1)$, depending on $n, \Omega,
C_0, K, \alpha$, such that
\beq\label{no-vorticity}
|u_\epsilon(x_0)|\ge 1-L\eta^\theta.
\eeq
\end{theorem}

\pf Throughout the proof, we assume that $\epsilon_0>0$ is chosen so that
\begin{equation}\label{small_epsilon0}
\max\big\{r_0(\Omega), \epsilon_0^\alpha\big\}\le\eta|\log\epsilon_0|.
\end{equation}
Since $0\le \alpha<1$, there exists
$\alpha_1>0$ such that $\alpha+\alpha_1<1$.
Let $x_0^*\in\partial\Omega$ such that $|x_0-x_0^*|=d_0\equiv{\rm{dist}}(x_0,\partial\Omega)$.
We divide the proof into two cases.\\
(a) $d_0\ge \frac12\epsilon^{\alpha+\alpha_1}$.  If $d_0\le \frac{r}3$,
then we have that
$$ B_{d_0}(x_0)\subset B_{2d_0}^+(x_0^*)\subset  B_{\frac{2r}3}^+(x_0^*)
\subset  B_r^+(x_0).$$
Hence, by the boundary monotonicity inequality \eqref{bdry_mono01.1}
in Theorem \ref{bdry_mono00} and \eqref{small_epsilon0}, we have that for $0<\epsilon\le \epsilon_0$
and $0<r\le r_0(\Omega)$,
\begin{eqnarray}\label{interior_small1}
\Phi_\epsilon(u_\epsilon; x_0, d_0)\equiv d_0^{2-n}\int_{ B_{d_0}(x_0)}e_\epsilon(u)
&\le& 2^{n-2}\Phi_\epsilon^+(u_\epsilon; x_0^*, 2d_0)\nonumber\\
&\le& 2^{n-2}e^{\Lambda r}\Phi_\epsilon^+(u_\epsilon; x_0^*, \frac{2r}3)+Cr\nonumber\\
&\le& C\Phi_\epsilon^+(u_\epsilon; x_0, r)+Cr_0(\Omega)\le C\eta|\log\epsilon|.
\end{eqnarray}

If $d_0>\frac{r}3$, then we have that $ B_{\frac{r}3}(x_0)\subset  B_r^+(x_0)$ and hence
\begin{equation}\label{interior_small2}
\Phi_\epsilon(u_\epsilon; x_0,\frac{r}3)\le 3^{n-2}\Phi_\epsilon^+(u_\epsilon; x_0, r)\le C\eta|\log\epsilon|.
\end{equation}
It is readily seen that \eqref{no-vorticity} follows from \eqref{interior_small1}, \eqref{interior_small2}, and
the interior $\eta$-compactness theorem by \cite{BBO2}.

\smallskip
\noindent(b) $d_0<\frac12\epsilon^{\alpha+\alpha_1}$. Since $r\ge \epsilon^\alpha$,
we have $d_0\le \frac{r}2$ so that $B^+_{\frac{r}2}(x_0^*)\subset  B_r^+(x_0)$.
Hence
\begin{equation}\label{bdry_small00}
\Phi_\epsilon^+(u_\epsilon; x_0^*, \frac{r}2)\le 2^{n-2} \Phi_\epsilon^+(u_\epsilon; x_0, r)\le 2^{n-2}\eta|\log\epsilon|.
\end{equation}
The proof, in the case (b), is divided into the following steps. For simplicity, we will assume
that $\Omega=\mathbb R^n_+$ throughout the presentation.

\medskip
\noindent{\bf Step 1}. {{\it Selection of good annuli}.} For $\delta\in (0, \frac1{16})$ to be chosen later, let $k\in\mathbb N$ be such
that
$$\epsilon^{\alpha+\alpha_1} {\delta}^{-(k+1)}\le \epsilon^\alpha,\ \  \epsilon^{\alpha+\alpha_1} {\delta}^{-(k+2)}\ge \epsilon^\alpha,$$
or, equivalently,
$$k+2\approx\Big[\alpha_1\frac{|\log\epsilon|}{|\log\delta|}\Big].$$
For $j=1,\cdots, k$, set
$$I_j=\big[\epsilon^{\alpha+\alpha_1}\delta^{-j+1}, \epsilon^{\alpha+\alpha_1}\delta^{-j}\big],$$

$$\displaystyle \mathcal C_j=B_{\epsilon^{\alpha+\alpha_1}{\delta}^{-j}}^+(x_0^*)\setminus
B_{\epsilon^{\alpha+\alpha_1}{\delta}^{-j+1}}^+(x_0^*)\ \ {\rm{and}}\ \ \mathcal D_j
={\Gamma_{\epsilon^{\alpha+\alpha_1}{\delta}^{-j}}(x_0^*)
\setminus\Gamma_{\epsilon^{\alpha+\alpha_1}\delta^{-j+1}}(x_0^*).}$$
To simplify the presentation, we write $u$ and $g$ for $u_\epsilon$ and $g_\epsilon$ respectively.
For all $0<\epsilon\le\epsilon_0$, it follows from (\ref{bdry_mono01.2})
and \eqref{bdry_mono1.2} that
\begin{eqnarray*}
&&\sum_{j=1}^{k-1}\int_{I_j}\tau^{1-n}\Big[\int_{B_\tau^+(x_0^*)}\frac{(1-|u|^2)^2}{2\epsilon^2}+\int_{\Gamma_\tau(x_0^*)}\frac{\lambda_\epsilon}2 |u-g|^2\Big]\,d\tau\\
&&\le e^{C\epsilon^\alpha} \Phi^+_\epsilon(u; x_0^*, \epsilon^{\alpha})
+C\epsilon^\alpha
\le 2\eta|\log\epsilon|,
\end{eqnarray*}
and
\begin{eqnarray*}
&&\sum_{j=1}^{k-1}\Big[\int_{\mathcal C_j} |x-x_0^*|^{2-n}\big(|\frac{\partial u}{\partial \nu}|^2+\frac{(1-|u|^2)^2}{2(n-2)\epsilon^2}\big)
+\int_{\mathcal D_j} |x-x_0^*|^{2-n}\frac{\lambda_\epsilon |u-g|^2}{2(n-2)}\Big]\\
&&\le e^{C \epsilon^{\alpha}}\widehat{\Phi}^+_\epsilon(u; x_0^*,\epsilon^{\alpha})+C\epsilon^{\alpha}
\le 2\eta |\log{\epsilon}|.
\end{eqnarray*}
Hence there exists $1\le j_0\le k-1$ such that
\beq\label{good_annual0}
\int_{I_{j_0}}\tau^{1-n}\Big[\int_{B_\tau^+(x_0^*)}\frac{(1-|u|^2)^2}{2\epsilon^2}+\int_{\Gamma_\tau(x_0^*)}\frac{\lambda_\epsilon}2 |u-g|^2\Big]\,d\tau
\le\frac{2\eta|\log\epsilon|}{k-1}
\le \frac{8\eta}{\alpha_1}|\log\delta|,
\eeq
and
\begin{eqnarray}\label{good_annual}
&&\int_{\mathcal C_{j_0}} |x-x_0^*|^{2-n}\big(|\frac{\partial u}{\partial \nu}|^2+\frac{(1-|u|^2)^2}{2(n-2)\epsilon^2}\big)
+\int_{\mathcal D_{j_0}} |x-x_0^*|^{2-n}\frac{\lambda_\epsilon |u-g|^2}{2(n-2)}\nonumber\\
&&\le \frac{2\eta|\log\epsilon|}{k-1}\le \frac{8\eta}{\alpha_1}|\log{\delta}|.
\end{eqnarray}
Set $$r_1=r_1(\epsilon):=\epsilon^{\alpha+\alpha_1}\delta^{-j_0} (\ge \epsilon^{\alpha+\alpha_1}).$$
Applying Fubini's theorem to (\ref{good_annual0}) and \eqref{good_annual}, we can find
$r_2\in [\frac{r_1}8, \frac{r_1}4]$  such that
\beq\label{good_radius0}
r_2^{2-n}\Big[\int_{B_{r_2}^+(x_0^*)}\frac{(1-|u|^2)^2}{2\epsilon^2}
+\int_{\Gamma_{r_2}(x_0^*)}\frac{\lambda_\epsilon}2 |u-g|^2\Big]
\le  \frac{64\eta}{\alpha_1}|\log\delta|,
\eeq
and
\begin{equation}\label{good_radius}
r_2^{3-n}\int_{S_{r_2}(x_0^*)}\big(|\frac{\partial u}{\partial |x-x_0^*|}|^2+\frac{(1-|u|^2)^2}{2(n-2)\epsilon^2}\big)
\le \frac{64\eta}{\alpha_1}|\log\delta|.
\end{equation}
Observe that by (\ref{bdry_mono1.1}) and \eqref{good_annual},  we also have
\begin{eqnarray}\label{bdry_mono2}
\widehat{\Phi}^+_\epsilon(u; x_0^*, r_1)&\le& \widehat{\Phi}^+_\epsilon(u; x_0^*, \delta r_1)+C\int_{\mathcal C_{j_0}} |x-x_0^*|^{2-n}\big(|\frac{\partial u}{\partial |x-x_0^*|}|^2+\frac{(1-|u|^2)^2}{2(n-2)\epsilon^2}\big)\nonumber\\
&&+C\int_{\mathcal D_{j_0}} |x-x_0^*|^{2-n}\frac{\lambda_\epsilon |u-g|^2}{2(n-2)}+C\epsilon^\alpha\nonumber\\
&\le& \widehat{\Phi}^+_\epsilon(u; x_0^*,\delta r_1)+\frac{64\eta}{\alpha_1}|\log\delta|+C\epsilon^{\alpha}.
\end{eqnarray}
\noindent{\bf Step 2}. {\it Estimation of the tangential energy of $u$ on $\displaystyle S_{r_2}(x_0^*)$}. Applying (\ref{tangential_radial2}),
with $r=r_2$, and (\ref{good_radius}), we have that
\begin{eqnarray}
\label{tangential_radial3}
r_2^{3-n}\int_{S_{r_2}(x_0^*)}|\nabla_T u|^2&\le& C\Big[r_2^{3-n}\int_{S_{r_2}(x_0^*)}\big|\frac{\partial u}{\partial |x-x_0^*|}\big|^2
+r_2^{2-n}\int_{B_{r_2}^+(x_0^*)}e_\epsilon(u)\Big]
+C\epsilon^\alpha\nonumber\\
&\le& \frac{64\eta}{\alpha_1}|\log\delta|+\widehat{\Phi}^+_\epsilon(u; x_0^*, r_1)
+C\epsilon^{\alpha}.
\end{eqnarray}

\noindent{\bf Step 3}. {\it Estimation of $\int_{B_{\delta r_1}^+(x_0^*)} e_\epsilon(u)$}.
The crucial part is to estimate
$\int_{B_{\delta r_1}^+(x_0^*)} |\nabla u|^2.$ Here we need to make some nontrivial
extensions of the ideas developed by \cite{BBO1, BBO2} on
(\ref{GL_WA1})$_1$ under the Dirichlet boundary condition to the weak anchoring condition (\ref{GL_WA1})$_2$.
Recall that $u$ solves
\begin{equation}\label{GZ2}
d^*(u\times du)=0 \ {\rm{in}}\ {B_{r_2}^+(x_0^*)}.
\end{equation}
Here $d$ denotes the exterior differential and $d^*=(-1)^{n+1}\star d\star$ denotes the co-exterior differential, and $\star$ denotes
the Hodge star operator.
Consider the equation for $\psi: B_{r_2}^+(x_0^*)\to \R^3$:
\begin{equation}\label{auxi-neuman}
\begin{cases}
\Delta\psi = 0 &\ {\rm{in}}\ B_{r_2}^+(x_0^*)\\
\frac{\partial\psi}{\partial\nu}=u\times \frac{\partial u}{\partial\nu} & \ {\rm{on}} \ S_{r_2}(x_0^*)\\
\frac{\partial\psi}{\partial\nu}=-u\times \lambda_\epsilon(u-g) & \ {\rm{on}}\ \Gamma_{r_2}(x_0^*).
\end{cases}
\end{equation}
By the weak anchoring condition (\ref{GL_WA1})$_2$, we have $\frac{\partial\psi}{\partial\nu}=
u\times\frac{\partial u}{\partial\nu}$ on $\partial B_{r_2}^+(x_0^*)$. Since (\ref{GZ2}) implies
$$\int_{\partial B_{r_2}^+(x_0^*)} u\times\frac{\partial u}{\partial\nu}=0,$$
we conclude that there exists a solution $\psi$ to the equation (\ref{auxi-neuman}).
From the standard elliptic theory, we have
\begin{eqnarray}\label{psi-estimate1}
&&r_2^{2-n}\int_{B_{r_2}^+(x_0^*)}|\nabla\psi|^2\le Cr_2^{3-n}\int_{\partial B_{r_2}^+(x_0^*)}\big|\frac{\partial\psi}{\partial\nu}\big|^2\nonumber\\
&&\le C\Big[r_2^{3-n}\int_{S_{r_2}(x_0^*)}|\frac{\partial u}{\partial |x-x_0^*|}|^2
+r_2^{3-n}\int_{\Gamma_{r_2}(x_0^*)}\lambda_\epsilon^2|u-g|^2\Big].
\end{eqnarray}
By using $r_2\lambda_\epsilon\le K$, (\ref{good_radius0}) and (\ref{good_radius}), this implies that
\begin{eqnarray}\label{psi-estimate2}
r_2^{2-n}\int_{B_{r_2}^+(x_0^*)}|\nabla\psi|^2&\le&
C\frac{\eta}{\alpha_1}|\log {\delta}|.
\end{eqnarray}
Let $\chi_{B_{r_2}^+(x_0^*)}$ denote the characteristic function of $B_{r_2}^+(x_0^*)$. Then it follows from
the definition of $\psi$ that we can verify
\begin{equation}\label{GZ3}
d^*\big[\chi_{B_{r_2}^+(x_0^*)}(u\times du-d\psi)\big]
=(-1)^{n+1}\star d\star \big[\chi_{B_{r_2}^+(x_0^*)}(u\times du-d\psi)\big]=0 \ \ {\rm{in}} \ \ \mathcal D'(\R^n_+).
\end{equation}
{Note that (\ref{GZ3}) is equivalent to stating that} $\star\big[\chi_{B_{r_2}^+(x_0^*)}(u\times du-d\psi)\big]$ is a closed $(n-1)$-form
on $\R^n_+$, i.e.,
\begin{equation}\label{GZ30}
d\Big(\star\big[\chi_{B_{r_2}^+(x_0^*)}(u\times du-d\psi)\big]\Big)=0 \ \ {\rm{in}} \ \ \mathcal D'(\R^n_+).
\end{equation}
Hence, by the Hodge decomposition theorem (see \cite{BBO2} Proposition A.8),
there exists  a co-closed $(n-2)$ form $\phi$ in $H^1(\R^n_+)\cap C^\infty\big(\R^n_+\setminus \overline{B_{r_2}^+(x_0^*)}\big)$ such that
\begin{equation}\label{GZ4}
\begin{cases}
\ \ d\phi=\star\big[\chi_{B_{r_2}^+(x_0^*)}(u\times du-d\psi)\big] & \ {\rm{in}}\ \R^n_+,\\
 \ d^*\phi=0 & \ {\rm{in}}\ \R^n_+,\\
 \big({\rm{i}}_{\partial \R^n_+}\big)^*\phi=0 & \ {\rm{on}}\ \partial\R^n_+,\\
|\phi(x)||x|^{n-1}\le C  &  \ {\rm{if}}\ |x| \ {\rm{is \ large}}.
\end{cases}
\end{equation}
Here ${\rm{i}}_{\partial\R^n_+}:\partial \R^n_+\to\R^n_+$ denotes the inclusion map.
Moreover, by (\ref{psi-estimate2}) it holds that
\begin{eqnarray}\label{l2-bound}
\big\|\nabla\phi\big\|_{L^2(\R^n_+)}^2&\le& C\Big[\big\|\nabla u\big\|_{L^2(B_{r_2}^+(x_0^*))}^2
+\big\|\nabla\psi\big\|_{L^2(B_{r_2}^+(x_0^*))}^2\Big]\nonumber\\
&\le& Cr_1^{n-2}\big[\widehat{\Phi}_\epsilon^+(u; x_0^*, r_1)+\frac{\eta}{\alpha_1}|\log{\delta}|\big].
\end{eqnarray}
For $0<\beta<\frac12$  to be chosen later, similar to \cite{BBO1, BBO2}, let
$f\in C^\infty\big(\mathbb R_+, \mathbb R_+)$ be such that
\begin{equation}\label{f-function}
f(t)=\begin{cases}\frac{1}{t} &\ {\rm{if}}\ t\ge 1-\beta,\\
1 & \ {\rm{if}}\ t\le 1-2\beta,\\
|f'(t)|\le 4 & \ {\rm{if}}\ 1-2\beta\le t\le 1-\beta.
\end{cases}
\end{equation}
Define $\widehat{f}:\R^n_+\to \big[1, \frac{1}{1-\beta}\big]$ by
\beq\label{de-vort}
\widehat{f}(x)=\begin{cases} f(|u(x)|) & \ {\rm{if}}\ x\in B_{r_2}^+(x_0^*),\\
1 & \ {\rm{if}}\ x\in \R^n_+\setminus B_{r_2}^+(x_0^*).
\end{cases}
\eeq
Since $\widehat{f}^2 u\times du= \widehat{f} u\times d(\widehat f u)$ in $ \mathbb B_{r_2}^+(x_0^*)$, we obtain that
\begin{eqnarray}\label{GZ5}
(-1)^{n+1}\Delta \phi &=& \star\big[d\big(\widehat{f}^2 \star d\phi\big)+d\big((1-\widehat{f}^2)\star d\phi\big)\big]\nonumber\\
&=& \star d\big( \chi_{\mathbb B_{r_2}^+(x_0^*)}\widehat{f}^2u\times du\big)-\star d(\chi_{\mathbb B_{r_2}^+(x_0^*)}
\widehat{f}^2 d\psi)+\star d\big((1-\widehat{f}^2)\star d\phi\big)\nonumber\\
&=&\star \big(\chi_{\mathbb B_{r_2}^+(x_0^*)} d(\widehat{f} u)\times d(\widehat{f} u)\big)
+\star\big(H^{n-1} \big|_{S_{r_2}(x_0^*)} \widehat{f}^2 u\times du\wedge d|x-x_0^*|\big)\nonumber\\
&&+\star d(\chi_{\mathbb B_{r_2}^+(x_0^*)}(1-
\widehat{f}^2) d\psi\big)+\star\big(H^{n-1}\big|_{S_{r_2}(x_0^*)}d|x-x_0^*|\wedge d\psi\big)\nonumber\\
&&+\star d\big((1-\widehat{f}^2)\star d\phi\big)\nonumber\\
&=&\omega_1+\omega_2+\omega_3+\omega_4+\omega_5\ \ {\rm{in}}\ \ \R^n_+.
\end{eqnarray}
Here $H^{n-1}\lfloor_{E}$ denotes $(n-1)$-dimensional Hausdorff measure restricted on a measurable
set $E\subset\R^n_+$.
As in \cite{BBO2} page 467, for $1\le i\le 5$, let $\phi_i\in H^1(\R^n_+, \Lambda^{n-2}(\R^n_+))$ be the solution of
\begin{equation}\label{phii-eqn}
\begin{cases}
\Delta\phi_i=\omega_i & \ {\rm{in}}\ \R^n_+,\\
\big({\rm{i}}_{\partial\R^n_+}\big)^*\phi_i=\big({\rm{i}}_{\partial\R^n_+}\big)^*(d^*\phi_i)=0 & \ {\rm{on}} \ \partial\R^n_+,\\
|\phi_i(x)|\rightarrow 0 & \ {\rm{as}}\ |x|\rightarrow +\infty.
\end{cases}
\end{equation}
Set $\Psi=\phi-\sum_{i=1}^5\phi_i$. Since $d^*\phi=0$ in $\R^n_+$,
we have that $\big({\rm{i}}_{\partial\R^n_+}\big)^*(d^*\phi) =0$ on $\partial\R^n_+$. Hence
\beq\label{difference}
\begin{cases}
\Delta\Psi=0 &\ {\rm{in}}\ \R^n_+,\\
\big({\rm{i}}_{\partial\R^n_+}\big)^*\Psi=\big({\rm{i}}_{\partial\R^n_+}\big)^*(d^*\Psi)=0 & \ {\rm{on}} \ \partial\R^n_+,\\
|\Psi(x)|\rightarrow 0 & \ {\rm{as}}\ |x|\rightarrow +\infty.
\end{cases}
\eeq
Hence by \cite{BBO2} Appendix, we conclude that $\Psi\equiv 0$ in $\R^n_+$, or equivalently,
\beq\label{split}
\displaystyle\phi=\sum_{i=1}^5\phi_i \ \ {\rm{in}}\ \ \R^n_+.
\eeq
It turns out that the estimation of the Dirichlet energies of $\phi_1, \phi_2, \phi_3,$ and $\phi_5$ can be
made, similar to that of \cite{BBO1, BBO2}, while the estimation of the Dirichlet energy of $\phi_4$
relies on the Rellich type  estimate of $\psi$, which is provided by Lemma \ref{rellich1} below.
For completeness, we sketch the estimates of $\phi_1,\phi_2, \phi_3, \phi_4, \phi_5$ as follows.

\medskip
\noindent (i) {\it Estimation of $\phi_1$}. We point out that this is the most difficult term to estimate.
Observe that for $x\in B_{r_1}^+(x_0^*)$,  we have
\beq\label{omega1-est0}
\big|d(\widehat{f}u)\times d(\widehat{f}u)(x)\big|\begin{cases}
=\displaystyle\big|d(\frac{u}{|u|})\times d(\frac{u}{|u|})(x)\big|=0 &
\ {\rm{if}}\  |u(x)|\ge 1-\beta,\\
\le \displaystyle\frac{C}{\epsilon^2}\le \frac{C(1-|u(x)|^2)^2}{\beta^2\epsilon^2}
& \ {\rm{if}} \ |u(x)|<1-\beta.
\end{cases}
\eeq
Thus, by (\ref{good_annual}), we have
\begin{equation}\label{omega1-est}
\big\|\omega_1\big\|_{L^1(\R^n_+)}=\big\|d(\widehat{f}u)\times d(\widehat{f}u)\big\|_{L^1(B_{r_2}^+(x_0^*))}
\le \frac{C}{\beta^2}\int_{B_{r_2}^+(x_0^*)} \frac{(1-|u|^2)^2}{\epsilon^2}.
\end{equation}
Now we need to estimate $\displaystyle\big\|\phi_1\big\|_{L^\infty(\R^n_+)}$. Similar to \cite{BBO2}, it follows
from \eqref{omega1-est0} that
for $x\in \R^n_+$,
\begin{eqnarray}\label{phi1-formula}
|\phi_1(x)|\le C\int_{\R^n_+}\frac{|\omega_1(y)|}{|x-y|^{n-2}}\,dy
\le \frac{C}{\beta^2}\int_{B_{r_2}^+(x_0^*)}|x-y|^{2-n}\frac{(1-|u(y)|^2)^2}{\epsilon^2}\,dy.
\end{eqnarray}
Since $\phi_1$ is harmonic in $\R^n_+\setminus \overline{ B_{r_2}^+}(x_0^*)$ and $|\phi_1(x)|\rightarrow 0$ as
$|x|\rightarrow \infty$, we have
\begin{equation}\label{max-phi1}
\big\|\phi_1\big\|_{L^\infty(\R^n)}=\max_{x\in\overline{ B_{r_2}^+(x_0^*)}}|\phi_1(x)|.
\end{equation}
The following Lemma plays a crucial role in the estimation of $\|\nabla\phi_1\|_{L^2(\mathbb R^n_+)}$.
\begin{lemma}\label{max-phi1-est1} Let $\phi_1$ be as above. There exists $C>0$, independent of $\epsilon$,
such that
\begin{equation}\label{max-phi1-est2}
\big\|\phi_1\big\|_{L^\infty(B_{r_2}^+(x_0^*))}\le \frac{C}{\beta^2}\big[\widehat{\Phi}^+_\epsilon(u; x_0^*, r_1)+\epsilon^{\alpha}\big].
\end{equation}
\end{lemma}
\noindent{\it Proof of Lemma \ref{max-phi1-est1}}:  For any $x=(x', x_n)\in  B_{r_2}^+(x_0^*)\setminus\Gamma_{r_2}(x_0^*)$, set $x^*=(x',0)\in \Gamma_{r_2}(x_0^*)$. From $x_n>0$, we have that
$$ B_{x_n}(x)\subset  B_{2x_n}^+(x^*)\subset  B_{\frac{r_1}2}^+(x^*)\subset  B_{r_1}^+(x_0^*).$$
Hence by Theorem 2.4, we obtain
\begin{eqnarray}\label{interior-mono1}
\widehat{\Phi}_\epsilon (u; x, x_n)\equiv x_n^{2-n}\int_{ B_{x_n}(x)}e_\epsilon(u)
&\le& 2^{n-2}\widehat{\Phi}_\epsilon^+(u; x^*, 2x_n)\nonumber\\
&\le& 2^{n-2}e^{C r_1}\big[\widehat{\Phi}_\epsilon^+(u; x^*, \frac{r_1}2)+r_1\big]\nonumber\\
&\le& C\big[\widehat{\Phi}_\epsilon^+(u; x_0^*, r_1)+\epsilon^{\alpha}\big].
\end{eqnarray}
This, combined with the interior monotonicity inequality in \cite{BBO2} Corollary II.1, gives
\begin{equation}\label{potential-estimate0}
\int_{B_{x_n}(x)}|x-y|^{2-n}\frac{(1-|u|^2)^2}{\epsilon^2}\,dy\le C\widehat{\Phi}_\epsilon (u; x, x_n)
\le C\big[\widehat{\Phi}_\epsilon^+(u; x_0^*, r_1)+\epsilon^{\alpha}\big].
\end{equation}
Observe that for any $y\in B_{r_2}^+(x_0^*)\setminus  B_{x_n}(x)$,
it follows from the triangle inequality and the condition $x_n\le |x-y|$ that $|x^*-y|\le 2|x-y|$.
Also observe that $B_{r_2}^+(x_0^*)\subset  B_{2r_2}^+(x^*)\subset  B_{\frac{r_1}2}^+(x^*)
\subset B_{r_1}^+(x_0^*)$, so that by Theorem 2.4 we have
\begin{equation}\label{bdry_mono55}\widehat{\Phi}_\epsilon^+(u; x^*, 2r_2)
\le C\big[\widehat{\Phi}_\epsilon(u; x^*, \frac{r_2}2)+\epsilon^{\alpha}\big]
\le C\big[\widehat{\Phi}_\epsilon(u; x_0^*, r_1)+\epsilon^{\alpha}\big].
\end{equation}
Hence by \eqref{phi1-formula}, \eqref{potential-estimate0}, \eqref{bdry_mono55}, and (\ref{bdry_mono5}), we obtain that
\begin{eqnarray*}|\phi_1(x)|&\le& \frac{C}{\beta^2}
\Big[\int_{ B_{x_n}(x)}|x-y|^{2-n}\frac{(1-|u(y)|^2)^2}{\epsilon^2}\,dy
+\int_{ B_{r_2}^+(x_0^*)\setminus B_{x_n}(x)}|x-y|^{2-n}\frac{(1-|u(y)|^2)^2}{\epsilon^2}\,dy\Big]\\
&\le& \frac{C}{\beta^2}
\Big[\int_{ B_{x_n}(x)}|x-y|^{2-n}\frac{(1-|u(y)|^2)^2}{\epsilon^2}\,dy
+\int_{B_{2r_2}^+(x^*)}|x^*-y|^{2-n}\frac{(1-|u(y)|^2)^2}{\epsilon^2}\,dy\Big]\\
&\le& \frac{C}{\beta^2}\big[\widehat{\Phi}^+_\epsilon(u; x_0^*, r_1)+\epsilon^{\alpha}\big].
\end{eqnarray*}
This, after taking the supremum over $x\in B_{r_2}^+(x_0^*)$, implies \eqref{max-phi1-est2}. Hence Lemma \ref{max-phi1-est1} is proved. \qed

\smallskip
It follows from \eqref{max-phi1-est2} and \eqref{max-phi1} that
\beq\label{bound-phi1}
\big\|\phi_1\big\|_{L^\infty(\R^n_+)}\le \frac{C}{\beta^2}\big[\widehat{\Phi}^+_\epsilon(u; x_0^*, r_1)+\epsilon^{\alpha}\big].
\eeq
Multiplying  the equation for $\phi_1$ by $\phi_1$, integrating the resulting equation over
$\R^n_+$,  and applying (\ref{omega1-est}) and (\ref{bound-phi1}), we obtain
\beq\label{energy-est-phi1}
\int_{\R^n_+}|\nabla\phi_1|^2\le \big\|\phi_1\big\|_{L^\infty(\R^n_+)}\big\|\omega_1\big\|_{L^1(\R^n_+)}
\le \frac{C}{\beta^4}\Big(\int_{B_{r_2}^+(x_0^*)} \frac{(1-|u|^2)^2}{\epsilon^2}\Big)
\Big[\widehat{\Phi}^+_\epsilon(u; x_0^*, r_1)+\epsilon^{\alpha}\Big].
\eeq

\medskip
\noindent (ii) {\it Estimation of $\phi_2$.} From \cite{BBO1} Appendix and \eqref{tangential_radial3}, we have that
\begin{eqnarray}\label{energy-est-phi2}
&&\int_{\R^n_+}|\nabla\phi_2|^2\le Cr_2\int_{S_{r_2}^+(x_0^*)}|\omega_2|^2\le Cr_2\int_{S_{r_2}^+(x_0^*)}|\nabla_T u|^2\nonumber\\
&&\le C\Big[r_1^{n-2}\alpha_1^{-1}\eta|\log{\delta}|+\epsilon^{\alpha} r_1^{n-2}
+r_1^{n-2}\widehat{\Phi}^+_\epsilon(u; x_0^*, r_1)\Big].
\end{eqnarray}
Since $\phi_2$ is harmonic in $\R^n_+\setminus S_{r_2}(x_0^*)$ and satisfies
$$\big({\rm{i}}_{\partial\R^n_+}\big)^*\phi_2=\big({\rm{i}}_{\partial\R^n_+}\big)^*(d^*\phi_2)=0  \ {\rm{on}} \ \partial\R^n_+,$$
it follows from (\ref{energy-est-phi2}) and the Harnack inequality  that
\begin{eqnarray}\label{phi2-est}
&&\int_{B_{\delta r_1}^+(x_0^*)}|\nabla\phi_2|^2\le C\delta^n\int_{B_{r_2}^+(x_0^*)}|\nabla\phi_2|^2\nonumber\\
&&\le C\Big[\delta^n r_1^{n-2}\alpha_1^{-1}\eta|\log{\delta}|+\epsilon^{\alpha}\delta^n r_1^{n-2}+\delta^n r_1^{n-2}
\widehat{\Phi}^+_\epsilon(u; x_0^*,r_1)\Big].
\end{eqnarray}

\medskip
\noindent (iii)  {\it Estimation of $\phi_3$.} Multiplying the equation for $\phi_3$ by $\phi_3$, integrating
the resulting equation over
$\R^n_+$,  and
applying integration by parts and (\ref{psi-estimate2}), we have
\begin{eqnarray}\label{est-phi3}
\int_{\R^n_+}|\nabla\phi_3|^2&\le& C\int_{B_{r_2}^+(x_0^*)}\big|1-\widehat{f}^2\big|^2|\nabla\psi|^2
\le C\beta^2 \int_{B_{r_2}^+(x_0^*)}|\nabla\psi|^2\nonumber\\
&\le& C\beta^2r_1^{n-2}\big(\eta|\log{\delta}|+\epsilon^{\alpha}\big).
\end{eqnarray}

\medskip
\noindent (iv)  {\it Estimation of $\phi_5$.} Multiplying the equation for $\phi_5$ by $\phi_5$,
integrating the resulting equation over $\R^n_+$, and
applying integration by parts and \eqref {l2-bound}, we obtain that
\begin{eqnarray}\label{est-phi5}
&&\int_{\R^n_+}|\nabla\phi_5|^2\le C\int_{ B_{r_2}^+(x_0^*)}\big|1-\widehat{f}^2\big|^2|\nabla\phi|^2\nonumber\\
&&\le C\beta^2 \int_{ B_{r_2}^+(x_0^*)}|\nabla\phi|^2\le C\beta^2r_1^{n-2}\big[\widehat{\Phi}_\epsilon^+(u;x_0^*, r_1)
+\eta|\log\delta|+\epsilon^{\alpha}\big].
\end{eqnarray}

\medskip
\noindent (v) {\it Estimation of $\phi_4$}.  Similar to (iii), we have
\beq\label{1est-phi4}
\int_{\R^n_+}|\nabla\phi_4|^2\le Cr_2\int_{S_{r_2}(x_0^*)}|\omega_4|^2
\le Cr_2\int_{S_{r_2}(x_0^*)}|\nabla_T\psi|^2.
\eeq
Since $\phi_4$ is harmonic in $\R^n_+\setminus S_{r_2}(x_0^*)$ and satisfies
$$\big({\rm{i}}_{\partial \R^n_+}\big)^*(\phi_4)=\big({\rm{i}}_{\partial \R^n_+}\big)^*(d^*\phi_4)=0
\ {\rm{on}}\ \partial \R^n_+,$$
we have, by the mean value inequality and (\ref{1est-phi4}), that
\beq\label{2est-phi4}
\int_{B_{\delta r_1}^+(x_0^*)}|\nabla\phi_4|^2\le C\delta^{n}\int_{B_{r_2}^+(x_0^*)}|\nabla\phi_4|^2
\le C\delta^n r_2\int_{S_{r_2}(x_0^*)}|\nabla_T\psi|^2.
\eeq
Applying the boundary condition of (\ref{auxi-neuman}) for $\psi$, (\ref{psi-estimate2}),
\eqref{good_radius}, and
Lemma \ref{rellich1} below, we obtain
\begin{eqnarray}
\label{3est-phi4}
&&\int_{B_{\delta r_1}^+(x_0^*)}|\nabla\phi_4|^2\le
C\delta^n r_2\int_{S_{r_2}(x_0^*)}|\nabla_T\psi|^2
\le C\delta^n \Big[r_2\int_{\partial B_{r_2}^+(x_0^*)}\big|\frac{\partial\psi}{\partial\nu}\big|^2
+\int_{ B_{r_2}^+(x_0^*)}|\nabla\psi|^2\Big]\nonumber\\
&&\le C\delta^nr_1^{n-2}\big(\eta|\log{\delta}|+\epsilon^{\alpha}\big)
+C\delta^n\Big[r_2\int_{S_{r_2}(x_0^*)}\big|\frac{\partial u}{\partial r}\big|^2+r_2\int_{\Gamma_{r_2}(x_0^*)}\lambda_\epsilon^2|u-g|^2\Big]\nonumber\\
&&\le C\delta^nr_1^{n-2}\big(\eta|\log\delta|+\epsilon^{\alpha}\big).
\end{eqnarray}
Substituting the estimates (\ref{energy-est-phi1}), (\ref{phi2-est}),  (\ref{est-phi3}), (\ref{est-phi5}), and (\ref{3est-phi4}) into (\ref{split}) and applying (\ref{good_radius0}) and (\ref{psi-estimate2}), we obtain
\begin{eqnarray}\label{decay_estimate}
&&(\delta r_1)^{2-n}\int_{B_{\delta r_1}^+(x_0^*)}|u\times \nabla u|^2\nonumber\\
&&\le(\delta r_1)^{2-n}\int_{ B_{\delta r_1}^+(x_0^*)}(|\nabla\phi|^2+|\nabla\psi|^2)\nonumber\\
&&\le C\big(\delta^2+\beta^2\delta^{2-n}\big)\big(\eta|\log\delta|+\epsilon^{\alpha}\big)
+C\beta^2\delta^{2-n}\big[\widehat{\Phi}_\epsilon^+(u; x_0^*, r_1)+\eta|\log{\delta}|+\epsilon^{\alpha}\big]\nonumber\\
&&\quad  +C\delta^2\Big[\eta|\log{\delta}|+\epsilon^{\alpha}+\delta \widehat{\Phi}^+_\epsilon(u; x_0^*, r_1)\Big]\nonumber\\
&&\quad +C\delta^{2-n}\beta^{-4}\Big(r_2^{2-n}\int_{B_{r_2}^+(x_0^*)}\frac{(1-|u|^2)^2}{\epsilon^2}\Big)
\Big[\widehat{\Phi}^+_\epsilon(u; x_0^*, r_1)+\epsilon^{\alpha}\Big]\nonumber\\
&&\le C\big(\delta^2+\beta^2\delta^{2-n}\big)\big(\eta|\log{\delta}|+\epsilon^{\alpha}\big)\nonumber\\
&&\quad +C\Big[\delta^3+\beta^2\delta^{2-n}+\delta^{2-n}\beta^{-4}\eta|\log{\delta}|\Big]
\widehat{\Phi}^+_\epsilon(u; x_0^*, r_1)\nonumber\\
&&\le C\Big[\delta^3+\beta^2\delta^{2-n}+\delta^{2-n}\beta^{-4}\eta|\log{\delta}|\Big]
\widehat{\Phi}^+_\epsilon(u; x_0^*, r_1)\nonumber\\
&&\quad +C\big(\delta^2+\beta^2\delta^{2-n}\big)\big(\eta|\log{\delta}|+\epsilon^{\alpha}\big).
\end{eqnarray}
Recall that
\beq\label{pt_id1}
4|u|^2|\nabla u|^2=4|u\times \nabla u|^2+\big|\nabla |u|^2\big|^2.
\eeq
Now we need

\medskip

\noindent{(vi)} {\it Estimation of $\displaystyle\int_{B_{r_2}^+(x_0^*)}\big|\nabla |u|^2\big|^2$}. Observe that
$|u|^2$ solves the equation
\begin{equation}\label{square_u_eqn}
-\Delta(1-|u|^2)+\frac{2(1-|u|^2)}{\epsilon^2} |u|^2=2|\nabla u|^2.
\end{equation}
Multiplying (\ref{square_u_eqn}) by $(1-|u|^2)$ and integrating the resulting equation over $B_{r_2}^+(x_0^*)$, we obtain
\begin{eqnarray}\label{square_u_est}
&&\int_{B_{r_2}^+(x_0^*)}\big|\nabla |u|^2\big|^2+\int_{B_{r_2}^+(x_0^*)}\frac{2(1-|u|^2)^2}{\epsilon^2}|u|^2\nonumber\\
&&=2\int_{B_{r_2}^+(x_0^*)}(1-|u|^2)|\nabla u|^2+\int_{\partial B_{r_2}^+(x_0^*)}(1-|u|^2)\frac{\partial |u|^2}{\partial \nu}.
\end{eqnarray}
It is not hard to estimate
\begin{eqnarray}\label{square_u_est1}
&&\int_{B_{r_2}^+(x_0^*)}(1-|u|^2)|\nabla u|^2\nonumber\\
&&=
\Big\{\int_{\{x\in  B_{r_2}^+(x_0^*): |u(x)|\ge 1-\beta^2\}}
+\int_{\{x\in  B_{r_2}^+(x_0^*): |u(x)|\le 1-\beta^2\}}\Big\}(1-|u|^2)|\nabla u|^2\nonumber\\
&&\le C\beta^2 \int_{B_{r_2}^+(x_0^*)}|\nabla u|^2+C\beta^{-2} \int_{B_{r_2}^+(x_0^*)}\frac{(1-|u|^2)^2}{\epsilon^2}.
\end{eqnarray}
Since $|u-g|\ge 1-|u|$ and $\displaystyle\frac{\partial u}{\partial \nu}=-\lambda_\epsilon(u-g)$
 on $\Gamma_{r_2}(x_0^*)$, we can apply (\ref{good_radius0}) and (\ref{good_radius}) to get
\begin{eqnarray}\label{square_u_est2}
&&\int_{\partial B_{r_2}^+(x_0^*)}(1-|u|^2)\frac{\partial |u|^2}{\partial \nu}
= 2\int_{S_{r_2}(x_0^*)}(1-|u|^2)|u|\big|\frac{\partial u}{\partial \nu}\big|
+2\int_{\Gamma_{r_2}(x_0^*)}(1-|u|^2)|u|\big|\frac{\partial u}{\partial \nu}\big|\nonumber\\
&&\le C\epsilon\Big(\int_{S_{r_2}(x_0^*)}\frac{(1-|u|^2)^2}{\epsilon^2}\Big)^\frac12\Big(\int_{S_{r_2}(x_0^*)}
\big|\frac{\partial u}{\partial \nu}\big|^2\Big)^\frac12+C\int_{\Gamma_{r_2}(x_0^*)}\lambda_\epsilon|u-g|^2\nonumber\\
&&\le C(1+\epsilon)r_2^{n-2}\big(\eta|\log{\delta}|+\epsilon^{\alpha}\big).
\end{eqnarray}
Putting (\ref{square_u_est1}) and (\ref{square_u_est2}) into (\ref{square_u_est}) and applying (\ref{good_radius0}), we obtain
\begin{eqnarray}
 \label{square_u_est3}
r_2^{2-n}\int_{B_{r_2}^+(x_0^*)}\big|\nabla |u|^2\big|^2
&\le& C\beta^2r_2^{2-n} \int_{B_{r_2}^+(x_0^*)}|\nabla u|^2+C\beta^{-2} r_2^{2-n} \int_{B_{r_2}^+(x_0^*)}\frac{(1-|u|^2)^2}{\epsilon^2}\nonumber\\
&&+C\big(\eta|\log{\delta}|+\epsilon^{\alpha}\big)\nonumber\\
&\le& C\beta^2\widehat{\Phi}_\epsilon^+(u; x_0^*, r_1)+C\beta^{-2}\big(\eta|\log{\delta}|+\epsilon^{\alpha}\big).
\end{eqnarray}
From (\ref{pt_id1}), we can combine (\ref{square_u_est3}) with (\ref{decay_estimate}) to obtain
\begin{eqnarray}
\label{decay_estimate1}
(\delta r_1)^{2-n}\int_{B_{\delta r_1}^+(x_0^*)}|u|^2|\nabla u|^2
&\le& C\Big[\delta^3+\beta^2\delta^{2-n}+\delta^{2-n}\beta^{-4}\eta|\log{\delta}|\Big]
\widehat{\Phi}^+_\epsilon(u; x_0^*, r_1)\nonumber\\
&+&C\big(\delta^2+\beta^{-2}\delta^{2-n}\big)\big(\eta|\log{\delta}|+\epsilon^{\alpha}\big).
\end{eqnarray}
On the other hand, we can estimate
\begin{eqnarray}
&&(\delta r_1)^{2-n}\int_{B_{\delta r_1}^+(x_0^*)}(1-|u|^2)|\nabla u|^2\nonumber\\
&&= (\delta r_1)^{2-n}\Big\{\int_{ B_{\delta r_1}^+(x_0^*)\cap\{|u|\le 1-\beta^2\}}
+\int_{ B_{\delta r_1}^+(x_0^*)\cap\{|u|\ge 1-\beta^2\}}\Big\}(1-|u|^2)|\nabla u|^2\nonumber\\
&&\le C\beta^2 (\delta r_1)^{2-n}\int_{ B_{\delta r_1}^+(x_0^*)}|\nabla u|^2
+C\beta^{-2}(\delta r_1)^{2-n}\int_{ B_{\delta r_1}^+(x_0^*)}\frac{(1-|u|^2)^2}{\epsilon^2}\nonumber\\
&&\le C\beta^2 (\delta r_1)^{2-n}\int_{ B_{\delta r_1}^+(x_0^*)}|\nabla u|^2
+C\beta^{-2}\delta^{2-n}\big(\eta|\log{\delta}|+\epsilon^{\alpha}\big), \label{decay_estimate2}
\end{eqnarray}
where we have used (\ref{good_radius0}) in the last step.

Adding (\ref{decay_estimate1}) and (\ref{decay_estimate2}), we obtain
\begin{eqnarray}\label{decay_estimate3}
(1-C\beta^2) (\delta r_1)^{2-n}\int_{B_{\delta r_1}^+(x_0^*)}|\nabla u|^2
&\le& C\Big[\delta^3+\beta^2\delta^{2-n}+\delta^{2-n}\beta^{-4}\eta|\log{\delta}|\Big]
\widehat{\Phi}^+_\epsilon(u; x_0^*, r_1)\nonumber\\
&&+C\big(\delta^2+\beta^{-2}\delta^{2-n}\big)\big(\eta|\log{\delta}|+\epsilon^{\alpha}\big).
\end{eqnarray}
This, combined with (\ref{good_radius0}) again, implies that
\begin{eqnarray}\label{decay_estimate4}
\widehat{\Phi}_\epsilon^+(u; x_0^*, \delta r_1)
&\le& C\Big[\delta^3+\beta^2\delta^{2-n}+\delta^{2-n}\beta^{-4}\eta|\log{\delta}|\Big]
\widehat{\Phi}^+_\epsilon(u; x_0^*, r_1)\nonumber\\
&&+C\big(\delta^2+\beta^{-2}\delta^{2-n}\big)\big(\eta|\log{\delta}|+\epsilon^{\alpha}\big),
\end{eqnarray}
provided $C\beta^2\le \frac12.$

Observe that by putting (\ref{decay_estimate4}) into (\ref{bdry_mono2}), we obtain
\begin{eqnarray}\label{decay_estimate5}
&&\Big[1-C\Big(\delta^3+\beta^2\delta^{2-n}+\delta^{2-n}\beta^{-4}\eta|\log{\delta}|\Big)\Big]
\widehat{\Phi}^+_\epsilon(u; x_0^*, r_1)\nonumber\\
&&\le C\big(\delta^2+\alpha_1^{-1}+\beta^{-2}\delta^{2-n}\big)(\eta|\log{\delta}|+\epsilon^{\alpha}\big).
\end{eqnarray}
Set $\eta_1=\big(\frac{1}{2C}\big)^{\frac{3n}{n+1}}$. For any $0<\eta\le \eta_1$,
we choose $\delta=\eta^{\frac{1}{3n}}$, $\beta=\eta^{\frac{n+1}{6n}}$,
 $\epsilon_2\equiv\big(\frac{\eta|\log\eta|}{n}\big)^{\frac{1}{\alpha}}$, and
 set $\epsilon_0$  to equal the {smallest constant among \eqref{small_epsilon0}
 and $\min\{\epsilon_1, \epsilon_2\}$, where} $\epsilon_1=\big(\frac{\eta}{2C}\big)^{\frac{2}{\alpha}}$. It is clear that (\ref{decay_estimate5}) yields
that for any $0<\epsilon\le\epsilon_0$,
\beq\label{decay_estimate6}
\Big[1-C\big(\eta^{\frac1{n}}+\eta^{\frac{n+1}{3n}}\log(\frac{1}{\eta})\big)\Big]
\widehat{\Phi}_\epsilon^+(u;x_0^*, r_1)
\le C\eta^{\frac{n+1}{3n}}\log(\frac{1}{\eta}).
\eeq
Set $\eta_2>0$ such that
$$C\big(\eta^{\frac1{n}}+\eta^{\frac{n+1}{3n}}\log(\frac{1}{\eta})\big)\le \frac12, \
\forall 0<\eta\le\eta_2.$$
Define $\eta_0>0=\min\{\eta_1, \eta_2\}$.
Then we have that for any $0<\eta\le\eta_0$, it holds
\beq\label{decay_estimate7}
\widehat{\Phi}_\epsilon^+(u; x_0^*, r_1)
\le C\eta^{\frac{n+1}{3n}}\log(\frac{1}{\eta}), \ \forall 0<\epsilon\le \epsilon_0.
\eeq
From \eqref{decay_estimate7},  we claim that
\begin{equation}\label{decay_estimate8}
\widehat{\Phi}_\epsilon^+(u; x_0, \epsilon)
\le C\eta^{\frac{n+1}{3n}}\log(\frac{1}{\eta}), \ \forall 0<\epsilon\le \epsilon_0.
\end{equation}
The proof of \eqref{decay_estimate8} can be divided into two cases:\\
(i) $d_0=|x_0-x_0^*|\le \epsilon$. Then we have
that $B_\epsilon^+(x_0)\subset  B_{2\epsilon}^+(x_0^*)$ so that
by Theorem 2.4 it holds
$$\widehat{\Phi}_\epsilon^+(u; x_0, \epsilon)
\le 2^{n-2}\widehat{\Phi}_\epsilon^+(u; x_0^*, 2\epsilon)
\le C\big[\widehat{\Phi}_\epsilon^+(u; x_0^*, r_1)+\epsilon^{\alpha}\big]
\le C\eta^{\frac{n+1}{3n}}\log(\frac{1}{\eta}).$$
This yields \eqref{decay_estimate8}. \\
(ii) $d_0=|x_0-x_0^*|>\epsilon$.  Since
$d_0\le \frac12\epsilon^{\alpha+\alpha_1}\le \frac{r_1}2$, we can apply Theorem 2.4 and
\cite{BBO2} Corollary II.1 to obtain
\begin{eqnarray*}
\widehat{\Phi}_\epsilon(u; x_0,\epsilon)&\le& \widehat{\Phi}_\epsilon(u; x_0,d_0)
\le 2^{n-2}\widehat{\Phi}_\epsilon^+(u; x_0^*,2d_0)\\
&\le& C\big[\widehat{\Phi}_\epsilon^+(u; x_0^*, r_1)+\epsilon^{\alpha}\big]
\le C\eta^{\frac{n+1}{3n}}\log(\frac{1}{\eta}).
\end{eqnarray*}
This also yields \eqref{decay_estimate8}.

We may assume $|u(x_0)|<1$. It follows from Lemma 3.1 that $|\nabla u(x)|\le \frac{C}{\epsilon}$ and hence
$$|u(x)-u(x_0)|\le \frac{C}{\epsilon}|x-x_0|\le \frac{1-|u(x_0)|}2,$$
provided $\displaystyle |x-x_0|\le \bar{r}:=\frac{\epsilon}{2C}(1-|u(x_0)|)(\le \epsilon).$

It is clear that
\begin{eqnarray*}
\int_{B_{\epsilon}^+(x_0)}(1-|u(x)|^2)^2&\ge&\int_{B_{\bar r}^+(x_0)}(1-|u(x)|^2)^2
\ge \int_{B_{\bar r}^+(x_0)}(1-|u(x)|)^2\\
&\ge& \big(\frac{1-|u(x_0)|}2\big)^2|B_{\bar r}^+(x_0)|
\ge \frac{1}{C}\epsilon^n (1-|u(x_0)|)^{n+2}.
\end{eqnarray*}
Therefore we have
$$(1-|u(x_0)|)^{n+2}\le C \epsilon^{2-n}\int_{ B_\epsilon^+(x_0)}\frac{(1-|u(x)|^2)^2}{\epsilon^2}
\le C \widehat{\Phi}_\epsilon^+(u; x_0, \epsilon)\le C\eta^{\frac{n+1}{3n}}\log(\frac{1}{\eta}).
$$
This implies that for any $\eta\le\eta_0$, there exist $L>0$ and $0<\theta<1$, independent of $\eta$,
such that
\begin{equation}\label{eta11}
1-|u(x_0)|\le L\eta^\theta.
\end{equation}
On the other hand, \eqref{eta11} automatically holds for $\eta\ge\eta_0$ by choosing
$L=\eta_0^{-1}>0$ and $\theta=1$. Thus the proof is now complete. \qed

\medskip
Now we show the estimate of $L^2$-tangential energy of $\psi$ on $\partial B_{r_2}^+(x_0^*)$. Namely,
\begin{lemma}\label{rellich1} There exists  $C=C(n)>0$ such that
\beq\label{rellich2}
r_2\int_{\partial  B_{r_2}^+(x_0^*)}\big|\nabla_T\psi\big|^2
\le C\Big[r_2\int_{\partial  B_{r_2}^+(x_0^*)}\big|\frac{\partial\psi}{\partial\nu}\big|^2
+\int_{ B_{r_2}^+(x_0^*)}|\nabla\psi|^2\Big].
\eeq
\end{lemma}
\pf
Since $B_{r_2}^+(x_0^*)$ is  strictly star-shaped with respect to some interior point $a_0$: there is
a constant $c(n)>0$ such that
\beq\label{star-shape}
(x-a_0)\cdot\nu(x)\ge c(n) r_2, \ \forall\ x\in \partial  B_{r_2}^+(x_0^*).
\eeq
Multiplying the equation of (\ref{auxi-neuman}) by $(x-a_0)\cdot\nabla\psi$, integrating the resulting equation
over $ B_{r_2}^+(x_0^*)$, and
applying integration by parts, we obtain
\begin{eqnarray*}
&&0=\int_{B_{r_2}^+(x_0^*)}(x-a_0)\cdot\nabla \psi\Delta\psi\\
&&\ =\int_{ B_{r_2}^+(x_0^*)}\nabla\cdot\big\langle(x-a_0)\cdot \nabla\psi,\nabla\psi\big\rangle
-|\nabla\psi|^2-(x-a_0)\cdot\nabla\big(\frac{|\nabla\psi|^2}2\big)\\
&&\ =\int_{\partial  B_{r_2}^+(x_0^*)}\big\langle(x-a)\cdot\nabla\psi, \frac{\partial\psi}{\partial\nu}\big\rangle
+\frac{n-2}{2}\int_{ B_{r_2}^+(x_0^*)}|\nabla\psi|^2-\int_{\partial  B_{r_2}^+(x_0^*)}(x-a_0)\cdot\nu(x)\frac{|\nabla\psi|^2}2.
\end{eqnarray*}
By H\"older's inequality we have
\begin{eqnarray*}
&&\int_{\partial  B_{r_2}^+(x_0^*)}\big\langle(x-a_0)\cdot\nabla\psi, \frac{\partial\psi}{\partial\nu}\big\rangle\\
&&=\int_{\partial  B_{r_2}^+(x_0^*)} (x-a_0)\cdot\nu(x)\big|\frac{\partial\psi}{\partial\nu}\big|^2
+\int_{\partial  B_{r_2}^+(x_0^*)}(x-a_0)\cdot T(x) \big\langle \nabla_T\psi, \frac{\partial\psi}{\partial\nu}\big\rangle\\
&&\le Cr_2 \int_{\partial  B_{r_2}^+(x_0^*)}\big|\frac{\partial\psi}{\partial\nu}\big|^2
+\frac{c(n)r_2}4\int_{\partial  B_{r_2}^+(x_0^*)}|\nabla_T\psi|^2.
\end{eqnarray*}
Thus by (\ref{star-shape}) we obtain
\begin{eqnarray*}
&&c(n)r_2\int_{\partial B_{r_2}^+(x_0^*)}\frac{|\nabla\psi|^2}2
\le\int_{\partial  B_{r_2}^+(x_0^*)}(x-a_0)\cdot\nu(x)\frac{|\nabla\psi|^2}2\\
&&=
\int_{\partial  B_{r_2}^+(x_0^*)}\big\langle(x-a_0)\cdot\nabla\psi, \frac{\partial\psi}{\partial\nu}\big\rangle
+\frac{n-2}{2}\int_{ B_{r_2}^+(x_0^*)}|\nabla\psi|^2\\
&&\le \frac{c(n)r_1}4\int_{\partial  B_{r_2}^+(x_0^*)}|\nabla_T\psi|^2
+Cr_2 \int_{\partial  B_{r_2}^+(x_0^*)}\big|\frac{\partial\psi}{\partial\nu}\big|^2
+\frac{n-2}{2}\int_{ B_{r_2}^+(x_0^*)}|\nabla\psi|^2.
\end{eqnarray*}
This clearly implies (\ref{rellich2}). \qed

\section{Estimate of the potential energy on approximate vortex sets}
\setcounter{equation}{0}
\setcounter{theorem}{0}

In this section, we will show that the potential energy over any approximate vortex set is uniformly bounded
for a solution to the equation \eqref{GL_WA1}.

For $0<\beta<\frac12$, define the closed subset $S_\beta^\epsilon\subset\overline\Omega$ by
$$S^\epsilon_\beta:=\Big\{x\in\overline\Omega: \ \big|u_\epsilon(x)\big|\le 1-\beta\Big\}.$$
Then we have
\begin{theorem} \label{potential_bound1}
There exists $C_\beta>0$ depending on $\Omega$,
$\beta$, $K, \alpha$, $C_0$, and $M$
such that if  $g_\epsilon\in C^2(\partial\Omega,\mathbb S^1)$ {satisfies the condition $({\bf G})$
and $u_\epsilon\in C^2(\overline\Omega,\R^2)$ is a solution to
\eqref{GL_WA1}, with $\lambda_\epsilon=K\epsilon^{-\alpha}$ for some
$K>0$ and $\alpha\in [0,1)$, satisfying
$$F_\epsilon(u_\epsilon,\Omega)\le M|\log\epsilon|,  \forall \epsilon \in (0,1]$$
for some $M>0$,} then
\beq\label{potential_bound2}
\int_{S_{\beta}^\epsilon}\frac{(1-|u_\epsilon|^2)^2}{\epsilon^2}\le C_\beta, \ \forall \epsilon\in (0,1].
\eeq
\end{theorem}
\pf
It follows from both the interior monotonicity inequality
(see \cite{BBO2} Lemma II.2) and the boundary monotonicity inequality
\eqref{bdry_mono01.2} that for any $\displaystyle x\in S_\beta^\epsilon$, we can find
$\displaystyle r_x\in\big(\epsilon^{\alpha},
\epsilon^{\frac{\alpha}2}\big)$ such that
$$
r_x^{2-n}\int_{B_{r_x}(x)\cap\Omega}\frac{(1-|u_\epsilon|^2)^2}{\epsilon^2}
\le C\frac{F_\epsilon(u_\epsilon,\Omega)+\epsilon^{\frac{\alpha}2}}{|\log\epsilon|}
\le CM\frac{|\log\epsilon|+\epsilon^{\frac{\alpha}2}}{|\log\epsilon|}\le CM.
$$
Applying the Besicovitch covering
theorem (cf. \cite{EG} and \cite{BBH}), there exists a countable  family of points
$\displaystyle\Sigma=\big\{x_i\big\}_{i=1}^\infty\subset S_\beta^\epsilon$
such that \\
(i) $S_\beta^\epsilon\subset \bigcup_{i=1}^\infty B_{r_{x_i}}(x_i)$, and\\
(ii) there exists a positive integer $N(n)$ such that
 $\displaystyle\big\{B_{r_{x_i}}(x_i)\big\}_{i=1}^\infty$ can be decomposed into
 $N(n)$-families $\mathcal B_k$ of disjoint balls
for $1\le k\le N(n)$.

It follows easily that
\beq\label{potential-est3}
\int_{S_{\beta}^\epsilon}\frac{(1-|u_\epsilon|^2)^2}{\epsilon^2}\le C\sum_{i=1}^\infty r_{x_i}^{n-2}.
\eeq
For any $x\in S_\beta^\epsilon$, let $\displaystyle\eta=\eta(x, \epsilon):=r_x^{2-n}\frac{F_\epsilon^+(u_\epsilon; x, r_x)}{|\log\epsilon|}$. Let $L>0$ and $\theta>0$ be the common constants given by Theorem \ref{eta-comp} and the interior $\eta$-compactness Theorem (cf. \cite{BBO2} Theorem 2).

Set $\eta_0=\big(\frac{\beta}{2L}\big)^{\frac{1}{\theta}}>0$. Then we have\\
\noindent{\it Claim 4.1}. {\it There exists $\epsilon_0>0$, depending on $\eta_0,\Omega, C_0, K, \alpha$,
such that for any $0<\epsilon\le\epsilon_0$ and $x\in S_\beta^\epsilon$, $\eta(x,\epsilon)>\eta_0$.}

Suppose that this claim were false. Then for any $\epsilon_0>0$, there exists $\epsilon\in (0, \epsilon_0]$
and $x_\epsilon\in S_\beta^\epsilon$ such that
$$r_x^{2-n}F_\epsilon^+(u_\epsilon; x_\epsilon, r_x)\le \eta_0 |\log\epsilon|.$$
Let $\epsilon_0>0$ be the constant given by Theorem \ref{eta-comp} that corresponds to $\eta=\eta_0$. Then Theorem \ref{eta-comp} implies that
$$|u_\epsilon(x_\epsilon)|\ge 1-L\eta_0^\theta=1-L\Big[\big(\frac{\beta}{2L})^{\frac{1}{\theta}}\Big]^\theta=1-\frac{\beta}2.$$
Hence $x_\epsilon\not\in S_\beta^\epsilon$. This contradicts the choice of $x_\epsilon$.

It follows from Claim 4.1 that there exists $\epsilon_0>0$ such that for any $0<\epsilon\le\epsilon_0$
and $x\in S_\beta^\epsilon$, we have
$$r_x^{2-n}F_\epsilon^+(u_\epsilon; x, r_x)\ge \big(\frac{\beta}{2L}\big)^{\frac{1}{\theta}}\big|\log{\epsilon}\big|.$$
For each family $\mathcal B_k$, $1\le k\le N(n)$, we have
\begin{eqnarray*}
&&\big(\frac{\beta}{2L}\big)^{\frac{1}{\theta}}\big|\log{\epsilon}\big|\sum_{B_{r_{x_i}}(x_i)\in\mathcal B_k} r_{x_i}^{n-2}
\le\sum_{B_{r_{x_i}}(x_i)\in \mathcal B_k} F_\epsilon^+(u_\epsilon; x_i, r_{x_i})\\
&&=\int_{\bigcup_{B_{r_{x_i}}(x_i)\in\mathcal B_k}B_{r_{x_i}}(x_i)}e_\epsilon(u_\epsilon)
+\int_{\big(\bigcup_{B_{r_{x_i}}(x_i)\in\mathcal B_k}B_{r_{x_i}}(x_i)\big)\cap\partial\Omega}\lambda_\epsilon|u_\epsilon-g|^2\\
&&\le F_\epsilon(u_\epsilon,\Omega)\le M\big|\log\epsilon\big|.
\end{eqnarray*}
Hence
$$
\sum_{B_{r_{x_i}}(x_i)\in\mathcal B_k} r_{x_i}^{n-2}\le M\big(\frac{2L}{\beta}\big)^{\frac{1}{\theta}}.
$$
Taking the sum over $1\le k\le N(n)$ yields
\beq\label{potential-est4}
\sum_{i=1}^\infty r_{x_i}^{n-2}\le C(n, L, \theta)\beta^{-\frac{1}{\theta}}.
\eeq
Putting (\ref{potential-est3}) and (\ref{potential-est4}) together implies that (\ref{potential_bound2})
holds for all $0<\epsilon\le\epsilon_0$.  Since $|u_\epsilon|\le 1$ in $\Omega$,
it is easy to see that \eqref{potential_bound2} automatically holds for $\epsilon\in [\epsilon_0,1]$.
Hence \eqref{potential_bound2} holds for all $\epsilon\in (0, 1]$. This completes the proof. \qed

\section{Global $W^{1,p}$-estimate of $u_\epsilon$}
\setcounter{equation}{0}
\setcounter{theorem}{0}

In this section, we will utilize the potential energy estimate, Theorem \ref{potential_bound1}, to show
the global $W^{1,p}$-bound of any solution to \eqref{GL_WA1} under the global energy bound
\eqref{global_bound}. To handle the weak anchoring boundary condition \eqref{GL_WA1}$_2$,
we apply the Hodge decomposition to the quantity $\displaystyle\ u_\epsilon\times \frac{du_\epsilon}{|du_\epsilon|^q}$
for some $0<q<1$,
rather than $u_\epsilon\times du_\epsilon$ by \cite{BBO2} in their proof of Theorem 1 (7) for the Dirichlet boundary condition.
It turns out that our approach also works for the Dirichlet boundary condition as well.

\begin{theorem}\label{Lp-estimate1} {Assume that $g_\epsilon\in C^2(\partial\Omega,\mathbb S^1)$ satisfies the condition
$({\bf G})$ and $u_\epsilon\in C^2(\overline\Omega,\R^2)$ is a solution to
\eqref{GL_WA1}, with $\lambda_\epsilon=K\epsilon^{-\alpha}$ for some
$K>0$ and $\alpha\in [0,1)$, satisfying
$$F_\epsilon(u_\epsilon,\Omega)\le M|\log\epsilon|, \forall \epsilon \in (0,1]$$
for some $M>0$. Then for any $1\le p<\frac{n}{n-1}$,
there exists $C_p>0$, independent of $\epsilon$,  such that
\beq\label{Lp-estimate2}
\int_\Omega |\nabla u_\epsilon|^p\le C_p, \ \forall \epsilon\in (0,1].
\eeq}
\end{theorem}
\pf
For any $\frac{n-2}{n-1}<q<1$, applying the Hodge decomposition theorem (cf. \cite{BBO2} Appendix) to
the $1$-form $\displaystyle u_\epsilon\times \frac{du_\epsilon}{|du_\epsilon|^q}$ in $\Omega$,
we conclude that there exist $F_\epsilon\in W^{1, \frac{2-q}{1-q}}_0(\Omega)$ and
$G_\epsilon\in W^{1, \frac{2-q}{1-q}}(\Omega, \Lambda^1(\R^n))$ such that
\beq\label{hodge1}
\begin{cases}
\displaystyle
u_\epsilon\times \frac{du_\epsilon}{|du_\epsilon|^q}=dF_\epsilon+d^*G_\epsilon, \ dG_\epsilon=0 & \ {\rm{in}}\ \Omega,\\
F_\epsilon=0, \ \big({\rm{i}}_{\partial\Omega}\big)^*G_\epsilon=0 & \ {\rm{on}}\ \partial\Omega,\\
\end{cases}
\eeq
and
\beq\label{hodge-est}
\big\|F_\epsilon\big\|_{W^{1,\frac{2-q}{1-q}}(\Omega)}+\big\|G_\epsilon\big\|_{W^{1,\frac{2-q}{1-q}}(\Omega)}
\le C(q)\big\|\nabla u_\epsilon\big\|_{L^{2-q}(\Omega)}^{1-q}.
\eeq

For $0<\beta<1$, let $f$ be given by (\ref{f-function}). Then we have
\begin{eqnarray}\label{splitting-est}
&&\int_\Omega {f}^2(|u_\epsilon|)\frac{|u_\epsilon\times du_\epsilon|^2}{|du_\epsilon|^{q}}\nonumber\\
&&=
\int_\Omega \big\langle{f}^2(|u_\epsilon|) u_\epsilon\times du_\epsilon,
u_\epsilon\times \frac{du_\epsilon}{|du_\epsilon|^q}\big\rangle\nonumber\\
&&=\int_\Omega \big\langle{f}^2(|u_\epsilon|) u_\epsilon\times du_\epsilon, dF_\epsilon
+d^*G_\epsilon\big\rangle\\
&&=I+II.\nonumber
\end{eqnarray}
Since $F_\epsilon=0$ on $\partial\Omega$, and
$$d^*(u_\epsilon\times du_\epsilon)=0\ \ {\rm{in}}\ \ \mathcal D'(\Omega),$$
it follows from integration by parts that
$$\int_\Omega \langle u_\epsilon\times du_\epsilon, dF_\epsilon\rangle=0.$$
Hence, by (\ref{hodge-est}), we have
\begin{eqnarray}\label{I-est}
I&=&\int_\Omega \big\langle\big({f}^2(|u_\epsilon|)-1\big) u_\epsilon\times du_\epsilon,
dF_\epsilon\big\rangle\nonumber\\
&\le& C\big\|{f}^2(|u_\epsilon|)-1\big\|_{L^\infty(\Omega)}
\big\|u_\epsilon\times du_\epsilon\big\|_{L^{2-q}(\Omega)}\big\|F_\epsilon\big\|_{W^{1, \frac{2-q}{1-q}}(\Omega)}\nonumber\\
&\le& C\beta\big\|u_\epsilon\times du_\epsilon\big\|_{L^{2-q}(\Omega)}^{2-q}.
\end{eqnarray}

To estimate $II$, we first observe, as in the proof of the $\eta$-compactness Theorem 3.2, that
\begin{eqnarray}\label{L1-est}
\Big|d\big[{f}^2(|u_\epsilon|) u_\epsilon\times du_\epsilon\big]\Big|
&\le&\begin{cases}0 & \ {\rm{if}}\ \big|u_\epsilon(x)\big|\ge 1-\beta\\
\displaystyle \frac{C}{\beta^2}\frac{\big(1-|u_\epsilon|^2\big)^2}{\epsilon^2}
& \ {\rm{if}}\ \big|u_\epsilon(x)\big|\le 1-\beta
\end{cases}\nonumber\\
&\le& C\chi_{\{|u_\epsilon(x)|\le 1-\beta\}} \frac{\big(1-|u_\epsilon|^2\big)^2}{\beta^2\epsilon^2}.
\end{eqnarray}
Choose $\displaystyle \frac{n-2}{n-1}<q<1$ so that $\displaystyle n<\frac{2-q}{1-q}$.
Hence, from (\ref{hodge-est}) and Sobolev's embedding theorem,
we have that
\beq\label{bound-G}
\big\|G_\epsilon\big\|_{L^\infty(\Omega)}
\le C\big\|G_\epsilon\big\|_{W^{1,\frac{2-q}{1-q}}(\Omega)}
\le C\big\|\nabla u_\epsilon\big\|_{L^{2-q}(\Omega)}^{1-q}.
\eeq
By integration by parts, (\ref{L1-est}) and (\ref{bound-G}) imply that
\begin{eqnarray}\label{II-est}
II&=&\int_\Omega \big\langle d\big[{f}^2(|u_\epsilon|) u_\epsilon\times du_\epsilon\big], G_\epsilon\big\rangle\nonumber\\
&\le &\Big\|d\big[{f}^2(|u_\epsilon|) u_\epsilon\times du_\epsilon\big]\Big\|_{L^1(\Omega)}
\big\|G_\epsilon\big\|_{L^\infty(\Omega)}\nonumber\\
&\le& \frac{C}{\beta^2}\Big[\int_{\{|u_\epsilon(x)|\le 1-\beta\}} \frac{\big(1-|u_\epsilon|^2\big)^2}{\epsilon^2}\Big]
\big\|G_\epsilon\big\|_{L^\infty(\Omega)}\nonumber\\
&\le& \frac{C}{\beta^2}\Big[\int_{\{|u_\epsilon(x)|\le 1-\beta\}} \frac{\big(1-|u_\epsilon|^2\big)^2}{\epsilon^2}\Big]
\Big\|\nabla u_\epsilon\Big\|_{L^{2-q}(\Omega)}^{1-q}\nonumber\\
&\le& C(\beta)\big\|\nabla u_\epsilon\big\|_{L^{2-q}(\Omega)}^{1-q},
\end{eqnarray}
where we have used (\ref{potential_bound2}) in the last inequality.

Putting the estimates (\ref{I-est}) and (\ref{II-est}) into (\ref{splitting-est}), we arrive at
\beq\label{splitting-est1}
\int_\Omega {f}^2(|u_\epsilon|)\frac{|u_\epsilon\times du_\epsilon|^2}{|du_\epsilon|^{q}}
\le C\beta\big\|u_\epsilon\times du_\epsilon\big\|_{L^{2-q}(\Omega)}^{2-q}
+C(\beta)\big\|\nabla u_\epsilon\big\|_{L^{2-q}(\Omega)}^{1-q}.
\eeq
Since $f\ge 1$, it follows that
$$\int_\Omega \big|du_\epsilon\big|^{2-q}\le\int_\Omega {f}^2(|u_\epsilon|)\big|du_\epsilon\big|^{2-q},$$
and that  by (\ref{pt_id1}) it holds
$$4\big|u_\epsilon\big|^2\big|du_\epsilon\big|^2
=4\big |u_\epsilon\times du_\epsilon\big|^2+\big|d|u_\epsilon|^2\big|^2\
\ {\rm{in}}\ \Omega.$$
Moreover, by Kato's inequality, we have
$$|du_\epsilon|\ge |d|u_\epsilon||\ \ {\rm{in}}\ \Omega.$$
Thus we obtain
\begin{eqnarray}\label{Lp-est1}
&&\int_\Omega \big|du_\epsilon\big|^{2-q}\nonumber\\
&&\le\int_\Omega {f}^2(|u_\epsilon|)|u_\epsilon|^2\big|du_\epsilon\big|^{2-q}
+\int_\Omega {f}^2(|u_\epsilon|)\big(1-|u_\epsilon|^2\big)\big|du_\epsilon\big|^{2-q}\nonumber\\
&&\le \int_\Omega {f}^2(|u_\epsilon|)\frac{|u_\epsilon\times du_\epsilon|^2}{|du_\epsilon|^{q}}
+\int_\Omega {f}^2(|u_\epsilon|)\big|d|u_\epsilon|\big|^{2-q}\nonumber\\
&&\quad+\int_\Omega {f}^2(|u_\epsilon|)\big(1-|u_\epsilon|^2\big)\big|du_\epsilon\big|^{2-q}\nonumber\\
&&\le C\beta\big\|u_\epsilon\times du_\epsilon\big\|_{L^{2-q}(\Omega)}^{2-q}
+C(\beta)\big\|\nabla u_\epsilon\big\|_{L^{2-q}(\Omega)}^{1-q}\nonumber\\
&&\quad+C\Big[
\int_\Omega \big|d|u_\epsilon|\big|^{2-q}+\int_\Omega {f}^2(|u_\epsilon|)
\big(1-|u_\epsilon|^2\big)\big|du_\epsilon\big|^{2-q}\Big].
\end{eqnarray}
Now we need to estimate the last two terms in the right hand side of (\ref{Lp-est1}). It is not hard to estimate
\begin{eqnarray}\label{Lp-est2}
&&\int_\Omega {f}^2(|u_\epsilon|)
\big(1-|u_\epsilon|^2\big)\big|du_\epsilon\big|^{2-q}\nonumber\\
&&\le \Big\{\int_{\{|u_\epsilon(x)|\le 1-\beta\}}
+\int_{\{|u_\epsilon(x)|\ge 1-\beta\}}\Big\} (1-|u_\epsilon|^2)|du_\epsilon|^{2-q}\nonumber\\
&&\le C\beta \int_\Omega |du_\epsilon|^{2-q}+C\int_{\{|u_\epsilon(x)|\le 1-\beta\}} (1-|u_\epsilon|^2)|du_\epsilon|^{2-q}\nonumber\\
&&\le C\beta \int_\Omega |du_\epsilon|^{2-q}+C\int_{\{|u_\epsilon(x)|\le 1-\beta\}} \frac{(1-|u_\epsilon|^2)}{\epsilon}
|du_\epsilon|^{1-q}\nonumber\\
&&\le C\beta \int_\Omega |du_\epsilon|^{2-q}+C\Big(\int_{\{|u_\epsilon(x)|\le 1-\beta\}} \frac{(1-|u_\epsilon|^2)^2}{\epsilon^2}\Big)^\frac12 \Big(\int_\Omega |du_\epsilon|^{2-2q}\Big)^\frac12\nonumber\\
&&\le C\beta \int_\Omega |du_\epsilon|^{2-q}+C\Big(\int_\Omega |du_\epsilon|^{2-2q}\Big)^\frac12\nonumber\\
&&\le C\beta \int_\Omega |du_\epsilon|^{2-q}+C\Big(\int_\Omega |du_\epsilon|^{2-q}\Big)^\frac{1-q}{2-q},
\end{eqnarray}
where we have used $ |du_\epsilon|\le \frac{C}{\epsilon}$,
H\"older's inequality, and (\ref{potential_bound2}) in the derivation above.

Utilizing the weak anchoring condition (\ref{GL_WA1})$_2$, we can adopt the argument by \cite{BBO2} Proposition VI.4
to estimate $\int_\Omega |d|u_\epsilon||^{2-q}$ as follows. First recall
that $\rho_\epsilon:=|u_\epsilon|$ solves
the equation (see, e.g., (\ref{square_u_eqn})):
\beq\label{GZ50}\displaystyle
-\Delta \rho_\epsilon^2+2|\nabla u_\epsilon|^2
=\frac{2}{\epsilon^2}\rho_\epsilon^2\big(1-\rho_\epsilon^2\big)\ \ {\rm{in}}\ \ \Omega.
\eeq
Denote the set
$$\Sigma_\epsilon=\Big\{x\in\Omega: \ \rho_\epsilon(x)\ge 1-\epsilon^{\alpha}\Big\}.
$$
Set
$$\overline{\rho}_\epsilon=\max\Big\{\rho_\epsilon, 1-\epsilon^{\alpha}\Big\}
$$
so that $\overline{\rho}_\epsilon=\rho_\epsilon$ on $\Sigma_\epsilon$ and
$$\displaystyle0\le 1-\overline{\rho}_\epsilon\le \epsilon^\alpha \ \ {\rm{in}}\ \ \Omega.
$$
Multiplying (\ref{GZ50}) by $\overline{\rho}_\epsilon^2-1$ and integrating over $\Omega$
yields
\begin{eqnarray}\label{test-rho1}
&&\int_\Omega \nabla\rho_\epsilon^2\nabla\overline{\rho}_\epsilon^2
+2\int_\Omega \frac{\rho_\epsilon^2(1-\rho_\epsilon^2)(1-\overline{\rho}_\epsilon^2)}{\epsilon^2}\nonumber\\
&&=2\int_\Omega |\nabla u_\epsilon|^2(1-\overline{\rho}_\epsilon^2)
+2\int_{\partial\Omega}\frac{\partial \rho_\epsilon^2}{\partial\nu}(1-\overline{\rho}_\epsilon^2).
\end{eqnarray}
From $\displaystyle\frac{\partial u_\epsilon}{\partial\nu}=-\lambda_\epsilon(u_\epsilon-g_\epsilon)$ on $\partial\Omega$,
we see that
$$\big|\frac{\partial\rho_\epsilon}{\partial\nu}\big|\le \lambda_\epsilon |u_\epsilon-g_\epsilon|
\ \ {\rm{on}}\ \ \partial\Omega.$$
Hence (\ref{test-rho1}) implies
\beq\label{test-rho2}
\int_{\Sigma_\epsilon}|\nabla\rho_\epsilon^2|^2\le
2\epsilon^\alpha\int_\Omega |\nabla u_\epsilon|^2+
2\epsilon^{\alpha}\int_{\partial\Omega}\lambda_\epsilon |u_\epsilon-g_\epsilon|.
\eeq
Since $\lambda_\epsilon=K\epsilon^{-\alpha}$ and
$$E_\epsilon(u_\epsilon)+\int_{\partial\Omega}\lambda_\epsilon|u_\epsilon-g|^2\le M_0|\log\epsilon|,$$
we have
$$\int_\Omega |\nabla u_\epsilon|^2\le M_0 |\log\epsilon|,$$
and
$$\int_{\partial\Omega}\lambda_\epsilon |u_\epsilon-g|\le
C\lambda_\epsilon^\frac12\big(\int_{\partial\Omega}\lambda_\epsilon |u_\epsilon-g|^2\big)^\frac12
\le C\epsilon^{-\frac{\alpha}2}\big(1+M_0|\log\epsilon|\big).
$$
Thus we obtain
\beq\label{test-rho3}
\int_{\Sigma_\epsilon}|\nabla\rho_\epsilon^2|^2\le C\epsilon^{\frac{\alpha}2}\big(1+M_0|\log\epsilon|\big).
\eeq
On the other hand, since $1-\rho_\epsilon^2\ge \epsilon^{\alpha}$ in $\Omega\setminus \Sigma_\epsilon$
and
$$\int_\Omega (1-\rho^2_\epsilon)^2\le M_0\epsilon^2 |\log\epsilon|,$$
it follows that
$$\big|\Omega\setminus \Sigma_\epsilon\big|\le M_0\epsilon^{2(1-\alpha)} |\log\epsilon|.$$
Hence by H\"older's inequality we have
\begin{eqnarray}\label{test-rho4}
\int_{\Omega\setminus \Sigma_\epsilon}|\nabla\rho_\epsilon^2|^{2-q}
&\le& \big(\int_\Omega |\nabla\rho_\epsilon|^2\big)^{\frac{2-q}2}\big|\Omega\setminus \Sigma_\epsilon\big|^{\frac{q}2}\nonumber\\
&\le& C|\log\epsilon|^{\frac{2-q}2} \big(\epsilon^{2(1-\alpha)} |\log\epsilon|\big)^{\frac{q}2}\nonumber\\
&\le& C\epsilon^{(1-\alpha)q}|\log\epsilon|.
\end{eqnarray}
Combining (\ref{test-rho3}) with (\ref{test-rho4}), we see that there exists $0<\mu<(1-\alpha)q$, depending only
on $\alpha, M_0,$ and $q$, such that for any $0<\epsilon<1$,
\beq\label{Lp-est3}
\int_\Omega |\nabla\rho_\epsilon^2|^{2-q}\le C\epsilon^\mu.
\eeq
Substituting (\ref{Lp-est2}) and (\ref{Lp-est3}) into (\ref{Lp-est1}), we obtain
\begin{eqnarray}\label{Lp-est4}
\int_\Omega \big|du_\epsilon\big|^{2-q}
\le C\epsilon^\mu+C\beta\int_\Omega \big|du_\epsilon\big|^{2-q}
+C\Big(\int_\Omega |du_\epsilon|^{2-q}\Big)^\frac{1-q}{2-q}.
\end{eqnarray}
This, combined with Young's inequality, yields
\beq\label{Lp-est5}
\int_\Omega \big|du_\epsilon\big|^{2-q}
\le C\epsilon^\mu+(C\beta+\frac12)\int_\Omega \big|du_\epsilon\big|^{2-q}+C(q).
\eeq
Therefore, by choosing $\beta$ in $(0, \frac12)$ sufficiently small, we get
$$\int_\Omega |\nabla u_\epsilon|^{2-q}\le C(q).$$
This yields (\ref{Lp-estimate2}), since $q>\frac{n-2}{n-1}$ implies that $p=2-q\in (1, \frac{n}{n-1})$.
The proof is complete.
\qed

\section{Proof of Theorem \ref{compactness1}}
\setcounter{equation}{0}
\setcounter{theorem}{0}
In this section, we will use Theorem 3.2, Theorem 4.1, and Theorem 5.1 to give a proof of the main Theorem
\ref{compactness1}.

With the help of Theorem 3.2, Theorem 4.1, and Theorem 5.1, we can argue, similar to  that by \cite{BBO2}, to
show all the conclusions of Theorem \ref{compactness1}, except (a2) and the boundary regularity for $\alpha=0$.
For the latter, we need the following Lemmas.

\begin{lemma}\label{eta-H1} {There exist $\epsilon_0>0$ and $C_4>0$, depending  only on
$\Omega$, $C_0$, $K$, and $M$, and $r_0=r_0(\Omega)>0$ such that for any $\eta>0$ and $K>0$,
if for all $0<\epsilon\le \epsilon_0$, we have $g_\epsilon\in C^2(\partial\Omega,\mathbb S^1)$ satisfies the condition $({\bf G})$,
$u_\epsilon\in C^2(\overline\Omega,\mathbb R^2)$ is a solution to
\begin{equation}\label{GL_WA3}
\begin{cases}
\Delta u_\epsilon+\frac{1}{\epsilon^2}(1-|u_\epsilon|^2)u_\epsilon= 0 &  {\rm{in}}\ \Omega,\\
\frac{\partial u_\epsilon}{\partial\nu}+K(u_\epsilon-g_\epsilon)=0 & {\rm{on}}\ \partial\Omega,
\end{cases}
\end{equation}
satisfying} \eqref{global_bound} with some $M>0$, and for $x_0\in\partial\Omega$,
\begin{equation}\label{small_energy1}
\Phi_\epsilon^+(u_\epsilon; x_0, r_0)\le \eta|\log\epsilon|,
\end{equation}
then
\begin{equation}\label{no-vortex}
|u_\epsilon|\ge \frac12 \ \ {\rm{in}}\ \ B_{\frac{r_0}2}^+(x_0),
\end{equation}
 and
\begin{equation}\label{H1-bound}
\Phi^+_\epsilon\big(u_\epsilon; x_0, \frac{r_0}4\big)\le C_4.
\end{equation}

\end{lemma}
\pf Similar to Theorem 3.2, \eqref{no-vortex} can be proved by combining the boundary $\eta$-compactness Theorem 3.2
and the interior $\eta$-compactness Theorem in \cite{BBO2}.
To prove \eqref{H1-bound}, observe that
by \eqref{no-vortex} we can write
$$u_\epsilon=\rho_\epsilon e^{i\phi_\epsilon}\ \ {\rm{in}}\ \  B_{\frac{r_0}2}^+(x_0),$$
where $\rho_\epsilon=|u_\epsilon|$ and $\phi_\epsilon$ is defined up to an integer
multiple of $2\pi$. For simplicity, we assume
\begin{equation}\label{mean_bound}
\frac{1}{| B_{\frac{r_0}2}^+(x_0)|}\int_{ B_{\frac{r_0}2}^+(x_0)} \phi_\epsilon\in [0,2\pi).
\end{equation}
Since $u_\epsilon\in C^\infty(\overline\Omega)$, we have that
$\rho_\epsilon, \phi_\epsilon \in C^\infty( B_{\frac{r_0}2}^+(x_0))$.
Applying Theorem 5.1, we have
\begin{equation}\label{W1p-bound6}
\int_{ B_{\frac{r_0}2}^+(x_0)}\big(|\nabla\rho_\epsilon|^p+|\nabla\phi_\epsilon|^p\big)\le C_p, \ \forall\ 1\le p<\frac{n}{n-1}.
\end{equation}
Moreover, $\phi_\epsilon$ solves
\begin{equation}\label{GL-WA6}
\begin{cases}
-{\rm{div}}(\rho_\epsilon^2\nabla\phi_\epsilon)=0 & {\rm{in}}\  B_{\frac{r_0}2}^+(x_0),\\
\frac{\partial\phi_\epsilon}{\partial\nu}+\frac{K}{\rho_\epsilon}{\rm{Im}}(g_\epsilon e^{-i\phi_\epsilon})=0  & {\rm{on}}\ \Gamma_{\frac{r_0}2}(x_0).
\end{cases}
\end{equation}
Since $\rho_\epsilon\ge \frac12$ in $ B_{\frac{r_0}2}^+(x_0)$ and $|\frac{K}{\rho_\epsilon}{\rm{Im}}(g_\epsilon e^{-i\phi_\epsilon})|\le 2K
$ on $\Gamma_{\frac{r_0}2}(x_0)$, we can apply the standard elliptic theory to obtain that
\begin{equation}\label{W12-phi}
\big\|\nabla\phi_\epsilon\|_{L^2( B_{\frac{3r_0}8}^+(x_0))}\le C_0.
\end{equation}
This can be achieved by the Moser iteration argument.
In fact, let $\xi\in C^\infty_0( B_{r_0}(x_0))$
be a cut-off function of $ B_{\frac{r_0}2}(x_0)$ and $(\phi_\epsilon)_{x_0,r_0}$ be the average of $\phi_\epsilon$ over
$ B_{r_0}^+(x_0)$. For $1<q<\frac{n}{n-1}$, multiplying \eqref{GL-WA6} by $\xi^2\big|\phi_\epsilon-(\phi_\epsilon)_{x_0,r_0}\big|^{q-2}\big(\phi_\epsilon-(\phi_\epsilon)_{x_0,r_0}\big)$,  integrating over
$ B_{r_0}^+(x_0)$, and applying integration by parts and Young's inequality, we can obtain
\begin{eqnarray*}
&&\int_{ B_{r_0}^+(x_0)}\big|\nabla(\xi|\phi_\epsilon-(\phi_\epsilon)_{x_0,r_0}\big|^{\frac{q}2})\big|^2\\
&&\le 4\int_{ B_{r_0}^+(x_0)}\rho_\epsilon^2\big|\nabla(\xi|\phi_\epsilon-(\phi_\epsilon)_{x_0,r_0}\big|^{\frac{q}2})\big|^2\\
&&\le C\int_{ B_{r_0}^+(x_0)}\rho_\epsilon^2|\nabla\xi|^2
|\phi_\epsilon-(\phi_\epsilon)_{x_0,r_0}|^q+K\int_{\Gamma_{r_0}(x_0)}\rho_\epsilon\xi^2|\phi_\epsilon-(\phi_\epsilon)_{x_0,r_0}|^{q-1}\\
&&\le Cr_0^{-2}\int_{ B_{r_0}^+(x_0)}|\phi_\epsilon-(\phi_\epsilon)_{x_0,r_0}|^q
+K\int_{\Gamma_{r_0}(x_0)}|\phi_\epsilon-(\phi_\epsilon)_{x_0,r_0}|^{q-1}\\
&&\le C(r_0,q),
\end{eqnarray*}
where we have used in the last step both \eqref{W1p-bound6} and the trace inequality
$$\int_{\Gamma_{r_0}(x_0)}|\phi_\epsilon-(\phi_\epsilon)_{x_0,r_0}|^{q-1}
\le C(r_0, q-1)\big\|\phi_\epsilon\big\|_{W^{1,q-1}( B_{r_0}^+(x_0))}\le C_q.$$
On the other hand, by Sobolev's embedding theorem we have that
$$
\Big(\int_{ B_{r_0}^+(x_0)}\xi^{2^*}|\phi_\epsilon-(\phi_\epsilon)_{x_0,r_0}|^{\frac{2^*q}2}\Big)^{\frac{2}{2^*}}
\le C\int_{ B_{r_0}^+(x_0)}\big|\nabla(\xi|\phi_\epsilon-(\phi_\epsilon)_{x_0,r_0}\big|^{\frac{q}2})\big|^2,
$$
where $2^*=\frac{2n}{n-2}$ is the Sobolev exponent. Therefore, we conclude that
$$
\int_{ B_{\frac{r_0}2}^+(x_0)}\big|\phi_\epsilon-(\phi_\epsilon)_{x_0,{\frac{r_0}2}}\big|^{\frac{2^*q}2}
\le C(r_0, q).$$
Repeating the same argument, with $q$ replaced by $\frac{2^{*}q}2$, on $ B_{\frac{r_0}2}^+(x_0)$
for finitely many times, we can finally obtain that
$$
\int_{ B_{\frac{7r_0}{16}}^+(x_0)}\big|\phi_\epsilon-(\phi_\epsilon)_{x_0,{\frac{r_0}2}}\big|^{2}
\le C(r_0,q).
$$
Now we can test the equation \eqref{GL-WA6} by $\xi^2(\phi_\epsilon-(\phi_\epsilon)_{x_0, \frac{7r_0}{16}})$,
with $\xi$ being a cut-off function of $B_{\frac{3r_0}8}(x_0)$, to obtain \eqref{W12-phi}.

As for $\rho_\epsilon$, we have
\begin{equation}\label{GL-WA7}
\begin{cases}
-\Delta\rho_\epsilon+\rho_\epsilon|\nabla\phi_\epsilon|^2=\frac{1}{\epsilon^2}(1-\rho_\epsilon^2)\rho_\epsilon & {\rm{in}}\  B_{\frac{r_0}2}^+(x_0),\\
\frac{\partial\rho_\epsilon}{\partial\nu}+K\big(\rho_\epsilon-{\rm{Re}}(g_\epsilon e^{-i\phi_\epsilon})\big)=0  & {\rm{on}}\ \Gamma_{\frac{r_0}2}(x_0).
\end{cases}
\end{equation}
Let $\xi$ be a cut-off function of $B_{\frac{r_0}4}(x_0)$. Multiplying the equation (\ref{GL-WA7}) by $(1-\rho_\epsilon)\xi$
and integrating over $ B_{\frac{r_0}2}^+(x_0)$, we obtain
\begin{eqnarray}\label{control-rho}
&&\int_{ B_{\frac{r_0}2}^+(x_0)}\big[|\nabla\rho_\epsilon|^2+\frac{\rho_\epsilon(1-\rho_\epsilon)^2}{\epsilon^2}(1+\rho_\epsilon)\big]\xi\\
&&\le\int_{ B_{\frac{r_0}2}^+(x_0)}\big[|\nabla\phi_\epsilon|^2\xi+|\nabla\rho_\epsilon||\nabla\xi|\big](1-\rho_\epsilon)
+K\int_{\Gamma_{\frac{r_0}2}(x_0)}|\rho_\epsilon-{\rm{Re}}(g_\epsilon e^{-i\phi_\epsilon})||1-\rho_\epsilon|\xi\nonumber\\
&&\le C_0,\nonumber
\end{eqnarray}
where we have used \eqref{W12-phi} and \eqref{W1p-bound6} in the last step. This yields
\begin{equation}\label{W12-rho}
\int_{ B_{\frac{r_0}4}^+(x_0)}\big(|\nabla\rho_\epsilon|^2+\frac{(1-\rho_\epsilon^2)^2}{\epsilon^2}\big)\le C_0.
\end{equation}
Putting (\ref{W12-phi}) and (\ref{W12-rho}) together yields (\ref{H1-bound}). This completes the proof. \qed

\begin{lemma}\label{H1-Holder} Under the same assumptions as in Lemma \ref{eta-H1},
there exist $\beta_0\in (0,\frac12)$ and $C_1>0$, depending on $\Omega, K, M, \ and \ \|g\|_{C^1(\partial\Omega)}$,
such that $u_\epsilon\in C^{\beta_0}( B_{\frac{r_0}{8}}^+(x_0))$ and
\begin{equation}\label{Holder-bound}
\big[u_\epsilon\big]_{C^{\beta_0}( B_{\frac{r_0}{16}}^+(x_0))}\le C_1.
\end{equation}
\end{lemma}
\pf From Lemma \ref{eta-H1}, we have that
$$\|u_\epsilon\|_{H^1( B_{\frac{r_0}4}^+(x_0))}\le C(x_0, r_0)<\infty,$$
and
$$\Phi_\epsilon^+(u_\epsilon; x, \frac{r_0}8)\le C(x_0, r_0), \ \forall\  x\in  B_{\frac{r_0}8}^+(x_0).$$
Therefore by Theorem \ref{eta-comp} we have
$$|u_\epsilon(x)|\ge 1-L\eta_\epsilon^\theta, \ \forall\ x\in B_{\frac{r_0}8}^+(x_0),
$$
where $$\displaystyle\eta_\epsilon\equiv\frac{C(x_0, r_0)}{|\log\epsilon|}.$$
In particular, $\rho_\epsilon\rightarrow 1$ uniformly in $ B_{\frac{r_0}8}^+(x_0)$ as $\epsilon\rightarrow 0$.

Now for $x\in  B_{\frac{r_0}8}^+(x_0)$ and $0<r\le \frac{r_0}8$, let $\psi_\epsilon\in H^1(B_r^+(x))$
solve
\begin{equation}\label{u-equation}
\begin{cases}\Delta \psi_\epsilon=0 & \ {\rm{in}}\  B_r^+(x),\\
\frac{\partial \psi_\epsilon}{\partial\nu}+\frac{K}{\rho_\epsilon}{\rm{Im}}(g_\epsilon e^{-i\phi_\epsilon})=0 &  \ {\rm{on}}\ \Gamma_r(x),\\
\psi_\epsilon=\phi_\epsilon & \ {\rm{on}}\ S_r(x).
\end{cases}
\end{equation}
By the boundary regularity theorem on elliptic equations with oblique boundary conditions (see
Lieberman-Trudinger \cite{GT}), we have that $\psi_\epsilon\in C^\beta(B_{\frac{3r}{4}}^+(x))$ for all
$\beta\in (0,1)$, and
\begin{equation}\label{morrey-decay1}
(\delta r)^{2-n}\int_{B_{\delta r}^+(x)}|\nabla\psi_\epsilon|^2
\le C\delta^{2\beta} r^{2-n}\int_{ B_r^+(x)}|\nabla \psi_\epsilon|^2,
\ \forall \ 0<\delta<\frac34.
\end{equation}
Now multiplying both equation \eqref{GL-WA6} and equation \eqref{u-equation} by $(\phi_\epsilon-\psi_\epsilon)$,
subtracting the resulting equations and
integrating over $B_r^+(x)$, we obtain
\begin{eqnarray*}
&&\int_{ B_r^+(x)}|\nabla(\phi_\epsilon-\psi_\epsilon)|^2\\
&&=\int_{ B_r^+(x)}(1-\rho_\epsilon^2)\nabla\phi_\epsilon\cdot \nabla(\phi_\epsilon-\psi_\epsilon)
+K\int_{\Gamma_r(x)}(1-\rho_\epsilon^2)\rho_\epsilon^{-1}{\rm{Im}}(g_\epsilon e^{-i\phi_\epsilon})(\phi_\epsilon-\psi_\epsilon)\\
&&\le L\eta_\epsilon^\theta \int_{ B_r^+(x)}|\nabla\phi_\epsilon|^2+C\eta_\epsilon^\theta r^{n+2}
+\frac12 \int_{ B_r^+(x)}|\nabla(\phi_\epsilon-\psi_\epsilon)|^2,
\end{eqnarray*}
where we have used the following inequality in the last step
\begin{eqnarray*}
&&\int_{\Gamma_r(x)}|\phi_\epsilon-\psi_\epsilon|
\le Cr\int_{ B_r^+(x)}|\nabla(\phi_\epsilon-\psi_\epsilon)|\le Cr^{\frac{n+2}2}
\Big(\int_{ B_r^+(x)}|\nabla(\phi_\epsilon-\psi_\epsilon)|^2\Big)^\frac12\\
&&\le C\eta_\epsilon^\theta r^{n+2}
+\frac12 \int_{ B_r^+(x)}|\nabla(\phi_\epsilon-\psi_\epsilon)|^2.
\end{eqnarray*}
Therefore we obtain
\begin{equation}\label{morrey-decay2}
\int_{ B_r^+(x)}|\nabla(\phi_\epsilon-\psi_\epsilon)|^2\le C\eta_\epsilon^\theta \int_{ B_r^+(x)}|\nabla\phi_\epsilon|^2+C\eta_\epsilon^\theta r^{n+2}.
\end{equation}
This, combined with \eqref{morrey-decay1}, implies
\begin{equation}\label{morrey-decay3}
(\delta r)^{2-n}\int_{ B_{\delta r}^+(x)}|\nabla\phi_\epsilon|^2\le
C\big(\delta^{2\beta}+\delta^{2-n}\eta_\epsilon^\theta\big) r^{2-n}\int_{B_r^+(x)}|\nabla\phi_\epsilon|^2+
C\delta^{2-n}\eta_\epsilon^\theta r^4.
\end{equation}
Returning to $\rho_\epsilon$, it is easy to see from \eqref{control-rho} that
\begin{equation}\label{control-rho1}
\int_{ B_{\frac{r}2}^+(x)}\big[|\nabla\rho_\epsilon|^2+\frac{(1-\rho_\epsilon^2)^2}{\epsilon^2}\big]
\le C(\eta_\epsilon^\theta+\epsilon^2)\int_{ B_r^+(x)}|\nabla  u_\epsilon|^2+C\eta_\epsilon^\theta r^{n-1}.
\end{equation}
Combining \eqref{morrey-decay3} with \eqref{control-rho1}, and choosing sufficiently small
$\delta$ and $\epsilon$, we obtain that
\begin{equation}\label{morrey-decay4}
(\delta r)^{2-n}\int_{ B_{\delta r}^+(x)}|\nabla u_\epsilon|^2\le \frac12
r^{2-n}\int_{ B_{r}^+(x)}|\nabla u_\epsilon|^2 +Cr
\end{equation}
holds for any $x\in  B_{\frac{r_0}8}^+(x_0)$ and $0\le r\le \frac{r_0}8$. It is well-known that
iterations of \eqref{morrey-decay4} yields that there exists $\beta_0\in (0, \frac12)$ such that
$$r^{2-n}\int_{ B_{r}^+(x)}|\nabla u_\epsilon|^2 \le Cr^{2\beta_0},
\ x\in  B_{\frac{r_0}8}^+(x_0), \ 0\le r\le \frac{r_0}8.$$
Hence by Morrey's decay lemma \cite{EG} we have that $u_\epsilon\in C^{\beta_0}( B_{\frac{r_0}{16}}^+(x_0))$
satisfies (\ref{Holder-bound}). The proof is complete.
\qed

\bigskip
\noindent{\bf Proof of Theorem \ref{compactness1}}.
It follows from Theorem \ref{Lp-estimate1} that for some $\epsilon_0\in (0,1]$,
\begin{equation}
\label{global-bound1}
\int_\Omega |\nabla u_\epsilon|^p \le C_p, \ \forall\ 1\le p<\frac{n}{n-1}
\end{equation}
holds for all $0<\epsilon\le \epsilon_0$.

Define a sequence of Radon measures $\mu_\epsilon$ on $\overline\Omega$ by
$$\mu_\epsilon:=\frac{1}{|\log\epsilon|}\Big(e_\epsilon(u_\epsilon)\,dx+\frac{\lambda_\epsilon}{2}|u_\epsilon-g_\epsilon|^2\,dH^{n-1}_{\partial\Omega}\Big),$$
where $dH^{n-1}_{\partial\Omega}$ denotes the $H^{n-1}$-measure restricted on $\partial\Omega$,
Then \eqref{global_bound} implies
\begin{equation}
\label{global-bound2}
\mu_\epsilon(\Omega)=\frac{F_\epsilon(u_\epsilon, \Omega)}{|\log\epsilon|}\le M.
\end{equation}
In view of \eqref{global-bound1} and \eqref{global-bound2}, for any sequence $\epsilon_i\rightarrow 0$,
we can extract a subsequence, still denoted as $\epsilon_i$, such that
\begin{equation}\label{1p-convergence}
u_{\epsilon_i}\rightharpoonup u_* \ {\rm{in}}\ W^{1,p}(\Omega), \ \forall\ 1\le p<\frac{n}{n-1},
\end{equation}
for some map $u_*\in W^{1,p}(\Omega,\mathbb S^1)$ for  all $1\le p<\frac{n}{n-1}$, and
\begin{equation}\label{measure-convergence}
\mu_{\epsilon_i}\rightharpoonup \mu_* \ {\rm{in}}\ (C(\overline\Omega))^*,
\end{equation}
for a non-negative Radon measure $\mu_*$ on $\overline\Omega$,
as weak convergence of Radon measures on $\overline\Omega$.

Multiplying the equation
\eqref{GL-WA0} by $\times u_{\epsilon_i}$, we see that
\begin{equation}\label{GL-WA10}
{\rm{div}}(\nabla u_{\epsilon_i}\times u_{\epsilon_i})=0 \ {\rm{in}}\ \Omega.
\end{equation}
Sending $\epsilon_i\rightarrow 0$, \eqref{GL-WA10}  and \eqref{1p-convergence} yield
\begin{equation}\label{GL-WA11}
{\rm{div}}(\nabla u_{*}\times u_{*})=0 \ {\rm{in}}\ \mathcal D'(\Omega).
\end{equation}
It follows from (\ref{global_bound}) that when $0<\alpha<1$, it holds
\begin{equation}\label{bdry-trace}
\int_{\partial\Omega}|u_\epsilon-g_\epsilon|^2\le \frac{2}{\lambda_{\epsilon_i}} M|\log\epsilon_i|
\le \frac{2M}{K}\epsilon_i^\alpha |\log\epsilon_i|\rightarrow 0, \ {\rm{as}}\ \epsilon_i\rightarrow 0.
\end{equation}
Since $u_\epsilon\rightarrow u_*$ in $L^2(\partial\Omega)$,  we then obtain that $\displaystyle\int_{\partial\Omega}|u_*-g_*|^2=0$ so that $u_*=g_*$ on $\partial\Omega$.

When $\alpha=0$, since $u_{\epsilon_i}$ satisfies
$$\int_\Omega \nabla u_{\epsilon_i}\times u_{\epsilon_i}\cdot\nabla\phi
+K\int_{\partial\Omega}(u_{\epsilon_i}-g_{\epsilon_i})\times u_{\epsilon_i}\phi=0, \ \forall\phi\in C^\infty(\overline\Omega),
$$
we have, after taking $i\rightarrow\infty$, that
\begin{equation}\label{distribution}
\int_\Omega \nabla u_*\times u_{*}\cdot\nabla\phi
-K\int_{\partial\Omega}g_*\times u_{*}\phi=0, \ \forall\phi\in C^\infty(\overline\Omega).
\end{equation}
This implies \eqref{GL-WA11} and $\big(\frac{\partial u_*}{\partial\nu}-Kg_*\big)\times u_*=0$ on $\partial\Omega$
in the distribution sense.
Hence we establish the part (a2).

It follows from the interior monotonicity inequality (see \cite{BBO2} Lemma II.2) and the boundary monotonicity inequality \eqref{bdry_mono01.1} for $u_\epsilon$  that
for any $x\in\overline\Omega$, $r^{2-n}\mu_*(B_r(x)\cap\Omega)$
is monotonically nondecreasing with respect to
$r>0$.
Hence
$$\Theta^{n-2}(\mu_*, x)=\lim_{r\rightarrow 0} r^{2-n}\mu_*(B_r(x)\cap\Omega)$$
exists for all $x\in\overline\Omega$ and is upper semicontinuous in $\overline\Omega$. Moreover,
from the $\eta$-compactness Theorem 1.2,
$$\Sigma=\Big\{x\in\overline\Omega: \Theta^{n-2}(\mu_*,x)\ge \eta\Big\}.
$$
It is standard (\cite{EG}) that $\Sigma$ is a closed set, with $H^{n-2}(\Sigma)<\infty$. Moreover, $\Sigma$ is ($n-2$)-rectifiable and
$\mu_* \ {\rm{L}}\ \Omega$ is a stationary varifold in $\Omega$. (See \cite{BBO2} Theorem IX.1.)

Also, for any $x_0\in \Omega\setminus\Sigma$, there exist $r_0>0$ and sufficiently large $i_0$ such that
$B_{r_0}(x_0)\subset\Omega\setminus\Sigma$ and
$$r_0^{2-n}\int_{B_{r_0}(x_0)}e_{\epsilon_i}(u_{\epsilon_i})\,dx
\le \eta{|\log\epsilon_i|},\ \forall \ i\ge i_0.$$
Thus by \cite{BBO2} Theorem VIII.2, $u_{\epsilon_i}\rightarrow u_*$ in $C^k(B_{\frac{r_0}2}(x_0))$ for
all $k\ge 1$. This implies that $u_*\in C^\infty(\Omega\setminus\Sigma,\mathbb S^1)$ is a smooth
harmonic map. Thus we prove (b1) for $0<\alpha<1$.

When $\alpha=0$, we have that for any $x_0\in\overline\Omega\setminus\Sigma$, there exists $r_0>0$ such that
$$r_0^{2-n}\Phi_{\epsilon_i}^+(u_{\epsilon_i}; x_0, r_0)\le \eta|\log\epsilon_i|, \ i\ge i_0.$$
We proceed as follows.\\
(1) ${\rm{dist}}(x_0, \Omega)\le \frac{r_0}2$. Then by Lemma \ref{eta-H1} and Lemma \ref{H1-Holder} there exist
$\theta_0\in (0,\frac12)$ and $C_0>0$ such that for all $i\ge i_0$,
$$[u_{\epsilon_i}]_{C^{\theta_0}(B_{\frac{r_0}2}^+(x_0))}\le C_0.$$
Therefore $u_{\epsilon_i}\rightarrow u_*$ in $C^{\theta_0}(B_{\frac{r_0}2}^+(x_0))$. \\
(2) ${\rm{dist}}(x_0, \Omega)>\frac{r_0}2$. Then
$$(\frac{r_0}2)^{2-n}\int_{B_{\frac{r_0}2}(x_0)}e_{\epsilon_i}(u_{\epsilon_i})=(\frac{r_0}2)^{2-n}\Phi_{\epsilon_i}^+(u_{\epsilon_i}; x_0, \frac{r_0}2)\le r_0^{2-n}\Phi_{\epsilon_i}^+(u_{\epsilon_i}; x_0, r_0)
\le \eta|\log\epsilon_i|, \ i\ge i_0.
$$
Then $u_{\epsilon_i}\rightarrow u_*$ in $C^k(B_{\frac{r_0}{4}}(x_0))$ for all $k\ge 1$, see \cite{BBO2} Theorem VIII.2.
Combining (1) with (2), we conclude that $u_{\epsilon_i}\rightarrow u_*$ in $C_{\rm{loc}}^k(\Omega\setminus\Sigma)
\cap C^{\theta_0}_{\rm{loc}}(\overline\Omega\setminus\Sigma)$ for any $k\ge 1$, and
$u_*\in C^\infty(\Omega\setminus\Sigma)\cap C^{\theta_0}(\overline\Omega\setminus\Sigma)$.
This completes the proof of Theorem \ref{compactness1}. \qed

\medskip
We close this paper by raising a few questions pertaining to Theorem \ref{compactness1}.

\begin{remark}{\rm
(a) Assume $0<\alpha<1$. What is the boundary regularity of $u_*$ on $\partial\Omega\setminus\Sigma$?
Is $u_*$ H\"older continuous up to $\partial\Omega\setminus\Sigma$?\\
(b) Assume $\alpha=0$. Does  $u_*$  have higher order regularity up to $\partial\Omega\setminus\Sigma$? In particular, is it true that $u_*\in C^2(\overline\Omega\setminus\Sigma,\mathbb S^1)$?}
\end{remark}

\noindent{\bf Acknowledgement}. {P. Bauman and D. Phillips are partially supported by NSF grant DMS-1412840. C. Wang is partially supported by NSF grant DMS 1522869.}


\end{document}